\documentclass[aos,preprint]{imsart}

\RequirePackage[OT1]{fontenc}
\RequirePackage{amsthm,amsmath,amssymb}
\RequirePackage[numbers]{natbib}
\RequirePackage[colorlinks,citecolor=blue,urlcolor=blue]{hyperref}
\RequirePackage{graphicx}
\RequirePackage{algorithm2e,algorithmic}
\RequirePackage{multirow,hhline}

\arxiv{arXiv:0000.0000}

\startlocaldefs
\numberwithin{equation}{section}
\theoremstyle{plain}
\newtheorem{theorem}{Theorem}[section]
\newtheorem{lemma}{Lemma}[section]

\newtheorem{definition}{Definition}[section]
\newtheorem{assumption}{Assumption}[section]
\newtheorem{remark}{Remark}[section]
\newtheorem{corollary}{Corollary}[section]

\endlocaldefs

\begin{document}

\begin{frontmatter}
\title{A Unifying Framework of  High-Dimensional Sparse Estimation with Difference-of-Convex (DC) Regularizations}
\runtitle{DC unified regularization framework}

\begin{aug}
\author{\fnms{Shanshan} \snm{Cao}\thanksref{t1}\ead[label=e1]{scao36@gatech.edu}},
\author{\fnms{Xiaoming} \snm{Huo}\thanksref{t1}\ead[label=e2]{huo@gatech.edu}}
\and
\author{\fnms{Jong-Shi} \snm{Pang}\thanksref{t2}
\ead[label=e3]{jongship@usc.edu}
\ead[label=u1,url]{http://www.foo.com}}

\thankstext{t1}{Georgia Institute of Technoloty}
\thankstext{t2}{University of South California}
\runauthor{S. Cao, et al.}

%
%
\end{aug}

\begin{abstract}
Under the linear regression framework, we study the variable selection problem when the underlying model is assumed to have a small number of nonzero coefficients (i.e., the underlying linear model is sparse).
Non-convex penalties in specific forms are well-studied in the literature for sparse estimation.
A recent work \cite{ahn2016difference} has pointed out that nearly all existing non-convex penalties can be represented as  difference-of-convex (DC) functions, which can be expressed as the difference of two convex functions, while itself may not be convex.
There is a large existing literature on the optimization problems when their objectives and/or constraints involve DC functions.
Efficient numerical solutions have been proposed.
Under the DC framework, directional-stationary (d-stationary) solutions are considered, and they are usually not unique.
In this paper, we show that under some mild conditions, a certain subset of d-stationary solutions in an optimization problem (with a DC objective) has some ideal statistical properties: namely, asymptotic estimation consistency, asymptotic model selection consistency, asymptotic efficiency.
The aforementioned properties are the ones that have been proven by many researchers for a range of proposed non-convex penalties in the sparse estimation.
Our assumptions are either weaker than or comparable with those conditions that have been adopted in other existing works.
This work shows that DC is a nice framework to offer a unified approach to these existing work where non-convex penalty is involved.
Our work bridges the communities of optimization and statistics.
\end{abstract}

\begin{keyword}[class=MSC]
\kwd[Primary ]{62J05}
\kwd{62F12}
\kwd[; secondary ]{62J12}
\end{keyword}

\begin{keyword}
\kwd{high-dimensional sparse regression}
\kwd{regularization}
\kwd{asymptotic optimality}
\kwd{LASSO estimator}
\end{keyword}

\end{frontmatter}

\section{Introduction}
\label{sec:intro}
Sparse estimation under a linear regression model is a fundamental and classical problem in statistics.
It continues to be highly active in the high-dimensional regime when the underlying parameter is believed to be sparse.
Properties on the resulting estimators have been extensively studied with different penalties of the sparsity in \cite{zhao2006model, wainwright2009sharp, fan2001variable, fan2004nonconcave, zou2006adaptive, huang2008adaptive, zhang2010nearly, zhang2010analysis, zhang2013multi}, etc. However, most existing works focus on the properties on a specific solution to the possibly nonconvex objective function, which is used to derive a sparse estimation of the unknown parameter. The stationary solutions  of other kind might also be of interest and possess satisfying properties, such as the desired asymptotic estimation consistency, asymptotic model selection consistency, asymptotic efficiency.
A unified framework for the penalized high-dimensional sparse estimation and the relation to a subfield of optimization problems, namely, the difference-of-convex (DC) programming are missing  in the literature.
We establish such a connection in this paper.

\subsection{Sparsity induced penalties}

We first present the formulation of high-dimensional sparse estimation in linear regression setting using sparsity induced penalties.
Consider observations $(y_1, x_1)$, $(y_2, x_2)$, $\ldots$, $(y_n, x_n)$, where we have the response $y_i \in \mathbb{R}$ and the predictor $x_i \in \mathbb{R}^p$ satisfy
\begin{equation*}
y_i = \beta^{\ast T}x_i + \epsilon_i.
\end{equation*}
Here, $\beta^{\ast T}$ is the transpose of the vector $\beta^{\ast} \in \mathbb{R}^p$, which is the true however unknown underlying parameter to be estimated.
We further assume that noises $\epsilon_i$'s are independently distributed, with $0$ mean and equal variance $ \sigma^2$ (which can be a sub-Gaussian distribution with variance parameter $\sigma^2$), and are independent of $x_i$'s.
The above model is commonly written in the following matrix form:
\begin{equation}
\label{eq:LR}
y = X\beta^{\ast} + \epsilon,
\end{equation}
where
the vector $y = (y_1, \cdots, y_n)^T \in \mathbb{R}^n$ is the response vector,
 $X \in \mathbb{R}^{n \times p}$ is the model matrix  with rows being individual predictors, $x_1^T$, $\cdots$, $x_n^T$,
and the random vector $\epsilon$ contains the noises.

In the high-dimensional regime where the number of the parameters, denoted by $p$, exceeds the sample size, denoted by $n$, one of the most important methods (according to many works such as \cite{candes2007dantzig, bickel2009simultaneous, buhlmann2011statistics}) is to estimate the parameter by using the LASSO \citep{tibshirani1996regression} approach.
It is interesting to note that a mathematically equivalent approach was proposed in \cite{chen1995examples} around the same time in the computational and applied mathematics literature.
LASSO is defined through solving the following convex optimization problem:
\begin{equation}
\label{eq:LASSO}
\hat{\beta}^{lasso}(y, X; \lambda) = \arg \min_{\beta \in \mathbb{R}^p}\Big\{\frac{1}{2n}\|y - X\beta \|_2^2 +\lambda \| \beta \|_1\Big\}
\end{equation}
The first term in the above objective is the goodness-of-fit measure (a.k.a., the residual sum of squares) in the linear regression model \eqref{eq:LR}.
The second term in the objective is a penalty function, which is the sum of absolute values: $ \sum_{i=1}^k \lambda |\beta_i|$.
We can further write the penalty in a more general form $\sum_{i=1}^k P_\lambda(\beta_i)$, where the univariate function $P_\lambda(x)$ takes the form $P_\lambda(x) = \lambda |x|$ in the LASSO approach.
Many existing works, including \cite{zhao2006model} and \cite{wainwright2009sharp} and others, have proved that with high probability (i.e., probability goes to $1$ as sample size goes to infinity), under some conditions on the design matrix and the choices for $\lambda$, the LASSO will be able to find the right signed support for the unknown parameter $\beta^{\ast}$.
The cases that have been studied include (1) when the matrix can be fixed or random, (2) the dimension of the unknown parameter is fixed or goes to infinity as the sample size increases, and (3) other interesting situations.

Despite the success of obtaining sparse estimators by using LASSO, it is also well known that the resulting estimator is biased.
This can be readily seen by considering the special scenario, where the design matrix $X$ is orthonormal, consequently the $L_1$ penalty leads to a soft-thresholding solution, which is biased from the true parameter $\beta^{\ast}$.
De-biasing procedure has been studied in \cite{zhang2014confidence, javanmard2014hypothesis, van2014asymptotically, javanmard2014confidence}.
In the present paper, we decide to focus on the regularization (i.e., adding a penalty function) approach, partially because the de-biasing approach may require solving multiple optimization problems, therefore could be computationally disadvantageous.
At the same time, we may explore the other algorithmic-design approaches in the future.

An effective extension of the  LASSO estimator is to replace the penalty function $P_\lambda(x)$  in \eqref{eq:LASSO} into some folded concave functions, which are non-convex.
Some representitive works include SCAD \cite{fan2001variable, fan2004nonconcave}, MCP \cite{zhang2010nearly}, adaptive LASSO \cite{zou2006adaptive}, capped-$l_1$ \cite{zhang2013multi}, together with others.
In general, this leads to an NP-hard problem; therefore no polynomial-time algorithm is known in finding the global optimal solution.
Specifically, SCAD is proposed in \citep{fan2001variable, fan2004nonconcave},  in order to debias the estimation when the parameter is numerically relatively large, which gives a constant penalty as the parameter is large enough.
Adaptive LASSO is studied in \citep{zou2006adaptive, huang2008adaptive}, which is motivated by the fact that in the orthogonal design, the bias of the parameter estimation is approximately $\lambda$ in LASSO. The authors suggest to give different sizes of penalties to different parameters, so that the variables with large coefficients have smaller weights in the $\ell_1$ penalty (depending on some consistent estimator $\hat{\beta}$ of $\beta^{\ast}$). Then they can reduce the estimation bias of the lasso, while retaining its sparsity property.
MCP is proposed in \citep{zhang2010nearly}, which also gives a constant penalty as the parameter is large enough.
Capped-$l_1$ in \citep{zhang2010analysis, zhang2013multi}  gives a penalty of truncated $l_1$ penalty to ensure a constant penalty when the estimation is large.
Consistency results, including measuring the squared distance of the estimation, prediction, signed support recovery, for the previously mentioned formulations can be found in the original works.

\subsection{Difference-of-convex (DC) unified penalties}

Recently, it is pointed out by \cite{ahn2016difference} that all  the previously listed penalties can be written in a unified DC form.
Especially, the first term is the $\ell_1$ penalty $\|\beta\|_1$.
This leads to a DC formulation; i.e., solving the penalized least square problem with a generalized DC penalty function.
In \cite{ahn2016difference},  they prove under some strong assumptions (strong convexity of the loss functions) that the d-stationary solutions found by a standard DC algorithm (i.e., DCA) is in fact the global minimal.
This result might not be surprising because under their assumptions, the objective function (the DC-penalized loss functions) in fact can be strongly convex, which makes the solution unique.
They also prove that the $\ell_0$ norm of the d-stationary solution has an upper-bound, which doesn't shed lights on the support recovery property.
Our paper is inspired by \cite{ahn2016difference}, compared to \cite{ahn2016difference}, we relax the assumptions on strong convexity for the loss function and prove  the existence of a class of d-stationary solutions, which have the oracle properties in the linear regression scenario.
Our result indicates that the assumption on the strong convexity of the loss function within the domain is not necessary.
In addition, the aforementioned work has an applied mathematics focus -- their $\ell_0$ norm result does not imply statistical properties of the d-stationary solutions.
In statistical literature, the distance between the d-stationary solution and the assumed ground truth is considered.
Our result will be more formulated towards the statistical properties of the d-stationary solutions: namely model estimation consistency, asymptotic convergence rate in estimation, and model selection consistency.
Despite the difference, it is interesting to point out that both work will require the restricted convexity assumption, which is assumed in nearly all related work.
Besides, we also generalize the results to DC penalized general loss functions.

The use of DC functions offers a general  framework on non-convex regularization.
Some special cases are discussed in \cite{loh2013regularized} and \cite{wang2014optimal}, although they don't explicitly mention the DC functions in their work.
The first work \cite{loh2013regularized} assumes that the penalty function $p_{\lambda}(\beta)$ is separable in parameter $\beta$ and each of the univariate penalty can be written as the difference of a convex function $p_{\lambda,\mu}(t)$ and a quadratic function $\frac{\mu}{2}t^2$, where $\mu$ is a known positive constant.
     Therefore one has $p_{\lambda}(t) = p_{\lambda,\mu}(t) - \frac{\mu}{2}t^2$.
They restrict the feasible region to a bounded region containing the ground truth $\beta^\ast$. Under certain regularity conditions on the penalty, including differentiability,  and restrictive strong convexity of the loss function, they give the optimal upper-bound for the estimation error as well as for prediction error.
Their assumption includes the popular studied penalties like SCAD and MCP.
On the other hand, They don't have results on the support recovery and they purposely eliminated possible stationary solutions outside the bounded feasible region they defined.
The second work \cite{wang2014optimal} mainly assumes the restricted strong convexity of the penalties and the loss functions.
They mainly discuss the elliptical design regression, least square loss, and logistic loss with SCAD, MCP penalties which can be written as the summation of the $\ell_1$ penalty $\|\beta\|_1$ and a concave function $q_{\lambda}(\beta)$ with proper bound in the concavity.
They argue that the local quadratic approximation algorithm they provided converges to a unique local minimum which enjoys the oracle properties as if you've already know the support for the true parameter.
They prove the estimation error upper-bound. And they are able to prove the support recovery results for linear regression model with least square loss function. Both are under the assumption that the concavity of the function $q_{\lambda}(\beta)$ is bounded.

Both works \cite{loh2013regularized, wang2014optimal} assume the decreasing first order derivative of the penalty function on the nonnegative real line, which is necessary for the unbiasedness for estimation of larger $\beta$.
They both restrict the penalties such that the objective function is strongly convex within some region where the local optimal solutions as well as the true unknown parameters are in the given convex set.

While in the current work, we solve the unconstrained problem and prove the asymptotic convergence results of the estimation for a class of local d-stationary solutions without using the assumption of the bounded concavity of the $q_{\lambda}(\cdot)$ function and constant penalty when the parameter is large enough, which allows us to include other penalties such as transformed $\ell_1$ \citep{zhang2014minimization, lv2009unified}  and logarithmic \cite{mazumder2011sparsenet}, into our analysis.
Equipped with the bounded convexity assumption, we further prove that the support recovery consistency for the class of d-stationary solutions we find near the ground truth.

From the computation perspective, there is a rich literature on solving the penalized (also known as regularized) problem numerically.
For example, Local Linear Approximation (LLA) in \cite{zou2008one} prove that one-step estimator from LLA performs well in SCAD with penalized likelihood estimation.
They also prove the asymptotic normality under some regularity conditions of the Fisher Information matrix.
In this paper, we would like to explore the relationship between LLA and the popular DC Algorithm (DCA) which is often used in DC programming.
It turns out that all the above mentioned algorithms are special cases of DCA.

This paper builds a bridge between optimization, where people focused on solving the optimization problem efficiently, and statistics, where people mainly focused on the inference (finding the estimation).
The link here would be the DC programming and DCA.
The DC programming enables us to generalize the classical penalized likelihood function to the DC penalized likelihood function, while DCA provides us efficient algorithms to solve the corresponding numerical problems.
Borrowing strength from existing literatures enables us to solve the optimization problem efficiently with convergence guarantees.
We unify the existing algorithms in the literature for finding the local minima of non-convex optimization problems under the DCA framework.

\subsection{Notations}
For a real number $q \in [1, +\infty)$, the $\ell_q$ norm of a vector $\beta = (\beta_1, \beta_2, \cdots, \beta_p) \in \mathbb{R}^p$ is defined as $\|\beta\|_q = (\sum_{i = 1}^p\beta_i^{q})^{1/q}$.
Specially, the $\ell_{\infty}$ norm is defined as $\|\beta\|_{\infty} = \max_{i = 1}^p \{|\beta_i|\}$.
The $\ell_{0}$ norm is defined as  $\|\beta\|_{0} = \mbox{card} \{\mbox{supp}(\beta)\}$,
where we have supp$(\beta) = \{i : \beta_i \neq 0\}$ and card$\{\cdot\}$ is the cardinality of the set.
We denote the cardinality of a set $S$ by $|S|$ and its complement by $S^c$.
For $\beta = (\beta_1, \beta_2, \cdots, \beta_p) \in \mathbb{R}^p$, we let $\beta_{S}$ denote the sub-vector (of $\beta$) whose elements correspond to the set $S$; we let $X_{S}$ denote the sub-matrix (of $X$) whose columns indices are correspond to the set $S$.

\subsection{Organization}

The rest of the paper is organized as follows.
We review basic properties of the DC programming and the DC functions in Section \ref{sec:dc-funcs}.
We form a penalized least square problem with a generalized DC-penalty in Section \ref{sec:Formulation}.
Under mild assumptions, we prove in Section \ref{MAINRESULTS} that a set of the d-stationary solutions are close to the ground truth.
Furthermore, they are also support recovery consistent.
We also extend our results to generalized loss functions, such as logistic loss, etc., in Section \ref{GLOSS}.
We provide the connections among popular exiting algorithms in DC programming and statistics estimation with non-convex objective functions in Section \ref{SEC:CCCP}.
We finally conclude this work in Section \ref{sec:conclusion}.
When possible, the technical proofs are relegated to the Appendix.

\section{DC functions and related basic properties}
\label{sec:dc-funcs}

In this section, we first provide the necessary background as well as a definition of the Difference-of-Convex (DC) functions, before proceeding to our formulation (Section \ref{sec:Formulation}) and the main results (Section \ref{MAINRESULTS} and \ref{GLOSS}).
We present the definition of DC functions and its known properties in Section \ref{sec:dc-properties}.
The directional derivatives  are reviewed in Section \ref{sec:d-derivative}.
The class of DC functions that we are particularly interested are reviewed in Section \ref{sec:relate-to-stat}.
We then define and study the directional stationarity (d-stationarity) that we focus on in this work in Section \ref{sec:d-stationary}.

\subsection{DC programming and DC functions}
\label{sec:dc-properties}

DC programming is pervasive nowadays in both optimization and statistics.
The DC program has been introduced in the literature from 1950's \citep{aleksandrov1950surfaces}.
The paper \cite{hartman1959functions} gives a wealth of basic properties of the DC functions, which are the functions that are used in the objectives and constraints in the DC programming.
In particular, the DC programming has been intensively studied in the field of optimization in the early 1980's \citep{horst1999dc, hiriart1985generalized,tao1997convex,tuy1987global}.
The following gives a formal definition for a DC function.
\begin{definition}
\label{def:dc-funct}
A function, $p(x)$, is called a \textit{difference-of-convex (DC)} function if we have
$$
p(x) = g(x) - h(x)
$$
where both $g(x)$ and $h(x) $  are convex functions.
\end{definition}
There are many known results regarding to DC functions and DC programming. We summarize these properties of DC programming from the literature in Appendix \ref{app:dcproperty}.

\subsection{Directional derivative}
\label{sec:d-derivative}

To enable our description, we define the {\it directional derivative} in the following.
\begin{definition}
For a function $Q(\beta)$ that is defined on $\Omega$ where $\beta \in \Omega \subset \mathbb{R}$ or $\mathbb{R}^p$, for $\beta_0$, $\beta_1 \in \Omega$, the directional derivative at $\beta_0$ in the direction of $\beta_1 - \beta_0$ is defined as follows:
$$
Q'(\beta_0; \beta_1 - \beta_0) = \lim_{\tau \rightarrow 0+}\frac{Q(\beta_0 + \tau(\beta_1 - \beta_0)) - Q(\beta_0)}{\tau},
$$
where $\tau \in \mathbb{R}_+$.
\end{definition}
To compute the directional derivative, when a function $P(\beta)$ is differentiable in $\mathbb{R}^p$, the directional derivative with regard to $\beta$ at $\beta_0$ is given below:
\[
P'(\beta_0; \beta - \beta_0) = \langle \nabla P(\beta_0), \beta - \beta_0 \rangle,
\]
where $\nabla P(\beta_0)$ is the gradient of the function $P(\beta)$ at $\beta_0$, and $\langle \cdot, \cdot \rangle$ represents the inner product.
When the function $P(\beta)$ is non-differentiable however convex in $\mathbb{R}^p$, given its sub-gradient set $\partial P(\beta_0)$ at $\beta_0$,  the directional derivative with regard to $\beta$ at $\beta_0$ can be written as \cite[Theorem 23.4]{rockafellar2015convex}:
\[
P'(\beta_0; \beta - \beta_0) = \max_{v \in \partial P(\beta_0)} \langle v, \beta - \beta_0 \rangle.
\]

The recent papers \cite{nouiehed2017pervasiveness} and \cite{ahn2016difference} discuss the pervasiveness of the existence of the DC functions as well as its relation to statistics.
Specifically, they have the following results considering the pervasiveness of DC functions.
\begin{lemma} \cite[Proposition 1]{nouiehed2017pervasiveness}
For any univariate continuous concave function $p$ that is defined on $\mathbb{R}_+$, the composite function $\theta(|t|) = p(|t|)$ is Difference-of-Convex (DC) on $\mathbb{R}$ if and only if $p'(0;+)$, the directional derivative of $p(t)$ at $0$, which can be written as follows,
$$
p'(0; +) = \lim_{\tau \rightarrow 0+}
\frac{p(\tau) - p(0)}{\tau},
$$
exists and is finite.
\end{lemma}
The above lemma leads to the realization that nearly all the folded-concave penalties \citep{fan2001variable, fan2004nonconcave, zhang2010nearly, zhang2010analysis, zhang2013multi} in the sparsity study nowadays belong to the DC family.
We articulate details in the subsequent subsection.

\subsection{Relation to statistics}
\label{sec:relate-to-stat}

Based on the definition of the DC functions (Definition \ref{def:dc-funct}), it has been realized (e.g.,  \cite{wang2013calibrating},\cite{wang2014optimal}, and \cite{ahn2016difference}) that many well-studied penalties, such as SCAD, MCP, capped $\ell_1$, transformed $\ell_1$, logarithmic, can be written as DC functions.
We articulate details in the following.
Throughout this paper, we consider penalties $P(\beta)$ that are separable in the (potentially multivariate) parameter $\beta$: $P(\beta) = \sum_{i = 1}^{p}p(\beta_i)$, where $\beta = (\beta_1,\ldots,\beta_p)\in\mathbb{R}^p$.
We argue that function $p(x)$ is a DC function for the popular existing penalties in the literature; that is, we have $p(x) = g(x) - h(x)$, where functions $g$ and $h$ are convex. More specifically, for the penalties that are of interests to us and are widely used in statistical inference, we always have $g(x) = |x|$ (or $g(x) = \lambda |x|$, when an algorithmic parameter $\lambda$ is involved), however the function $h(x)$ varies per different penalties.

\begin{figure}[t]
  \begin{center}
   \includegraphics[width=.4\linewidth, height = .34\linewidth]{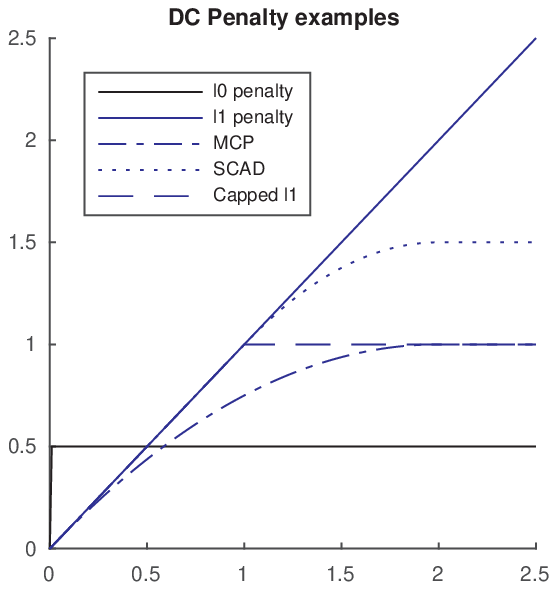}
   \includegraphics[width=.4\linewidth, height = .34\linewidth]{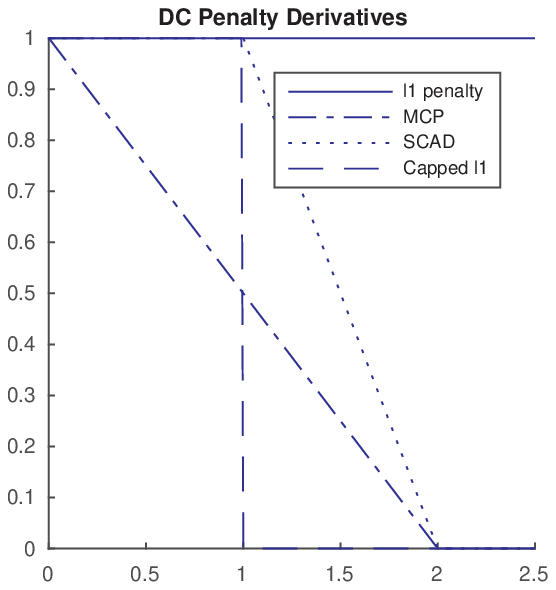}
  \caption{Examples of famous DC penalties and their derivatives in the literature: the $\ell_0$ penalty is plotted in the solid black line; the $\ell_1$ penalty function is plotted in the solid blue line; the MCP penalty is plotted in the dotted blue line; the MCP penalty is plotted in the dash-dot line; the Capped-$\ell_1$ penalty is plotted in the dashed blue line.}
    \label{FIG:PENALTY}
  \end{center}
\end{figure}

In the following, we describe how the popular penalty functions $p(x)$ mentioned previously can be decomposed as DC functions.
For simplicity, without loss of generality, we set the tuning parameter as $\lambda = 1$.
\begin{enumerate}
\item In {SCAD \citep{fan2001variable, fan2004nonconcave}}, we have
\[
p^{SCAD}(t) = |t| - h^{SCAD}_{\gamma}(t),
\]
where
\[
h^{SCAD}_{\gamma}(t) = \begin{cases}
      0 & |t| \leq 1 \\
      \frac{(|t|-1)^2}{2(\gamma-1)} & 1 \leq |t| \leq \gamma \\
      |t| - \frac{(\gamma+1)}{2} & |t| \geq \gamma
   \end{cases}
\]
and the function $h^{SCAD}_{\gamma}(t)$ can be verified to be convex on the positive real line and have a continuous first order derivative.

\item
In {MCP \citep{zhang2010nearly}}, we have
\[
p^{MCP}(t) = |t| - h^{MCP}_{\gamma}(t),
\]
where
\[
h^{MCP}_{\gamma}(t) = \begin{cases}
      \frac{|t|^2}{2\gamma} & |t| \leq \gamma \\
      |t| - \frac{\gamma}{2} & |t| \geq \gamma
   \end{cases}
\]
and the function $h^{MCP}_{\gamma}(t)$ can be verified to be convex on the positive real line and have a continuous first order derivative.

\item
In {Capped $\ell_1$ \citep{zhang2010analysis, zhang2013multi}}, we have
\[
p^{\text{capped } \ell_1}(t) = |t| - \max\left\{0,  \frac{2t}{\gamma} - 1, -\frac{2t}{\gamma} - 1\right\},
\]
where one can verify that both $|t|$ and $\max\left\{0,  \frac{2t}{\gamma} - 1, -\frac{2t}{\gamma} - 1\right\}$ are convex functions of $t$.

\item
In {Transformed $\ell_1$ \citep{zhang2014minimization, lv2009unified}}, we have
\[
p^{\text{Transformed }\ell_1}(t) = \frac{a+1}{a}|t| - \left(\frac{a+1}{a}|t| - \frac{(a+1)|t|}{a+ |t|}\right),
\]
where one can verify that both $\frac{a+1}{a}|t|$ and $\left(\frac{a+1}{a}|t| - \frac{(a+1)|t|}{a+ |t|}\right)$ are convex functions.

\item
In {Logarithmic \cite{mazumder2011sparsenet}}, we have
\[
p^{Log}(t) =  \frac{1}{\epsilon}|t| - \left(\frac{|t|}{\epsilon} - \log(|t|+\epsilon) + \log\epsilon\right),
\]
where $\epsilon$ is a given positive scalar; similarly, one can verify that both $\frac{\lambda}{\epsilon}|t|$  and $\left(\frac{|t|}{\epsilon} - \log(|t|+\epsilon) + \log\epsilon\right)$ are convex functions.
\end{enumerate}
Table \ref{table:penalty} summarizes the aforementioned DC decompositions. The  penalty functions and their first order derivatives in terms of $|t|$ are plotted in Figure  \ref{FIG:PENALTY}.
One common point that most of the above penalties share is that their first order derivative goes to 0 as $|t| \to \infty$.
\begin{table}[htbp]
\begin{center}
\begin{tabular}{|c|c|c|c|}
\hline
Penalty & $p(t)$ & $g(t)$  & $h(t)$\\
\hline
$\ell_1$ & $|t|$ & $|t|$ &  0  \\
\hline
SCAD&$\int^{|t|}_{0}1\wedge(1-\frac{x-1}{\gamma-1})_+dx$ & $|t|$ & $\frac{(|t|-1)^2}{2(\gamma-1)}I\{1<|t|<\gamma\} $   \\
 & &&$+ (|t| - \frac{\gamma+1}{2})I\{|t|\geq\gamma\}$ \\
 \hline
MCP&$\int^{|t|}_{0}(1-\frac{x}{\gamma})_+dx$& $|t|$& $(|t|- \gamma/2)I\{|t|>\gamma\}$  \\
 & & &$+\frac{t^2}{2\gamma}I\{|t|<\gamma\}$   \\
 \hline
Capped-$\ell_1$&$\min\{\gamma/2,|t|\}$& $|t|$& $\max\{0,\frac{2t}{\gamma}-1\}$   \\
\hline
Transformed $\ell_1$&$\frac{(a+1)|t|}{a+ |t|}$& $\frac{a+1}{a}|t|  $& $\frac{a+1}{a}|t| - \frac{(a+1)|t|}{a+ |t|} $  \\
\hline
Logarithmic&$\log(|t|+\epsilon) - \log\epsilon$& $\frac{|t|}{\epsilon} $ & $\frac{|t|}{\epsilon} - \log(|t|+\epsilon) + \log\epsilon$  \\
\hline
\end{tabular}
\caption{The DC decompositions of some well-known penalty functions in statistical inference.
The first column contains the name of the methodology.
The second column describes the penalty function.
The last two columns present the corresponding two convex functional components (i.e., $g$ and $h$) in the DC decomposition:
$p(t) =g(t) -h(t)$. }
\label{table:penalty}
\end{center}
\end{table}

Although there are many other different DC decompositions, this one has the advantage of easy interpretation and correcting the penalty of LASSO, which in some sense, does a debiasing job for LASSO estimator (by choosing $h_{\lambda}(t)$ to be linear with slope $\lambda$ when $|t|$ is large enough).
It also shares common features with popular penalties in literature, like SCAD, MCP, capped $\ell_1$ penalties where the penalty is close to or equal to $\ell_1$ penalty when the solution is around the origin.
Furthermore, the resulting penalty $p(x) = g(x) - h(x)$ is singular at $0$, which makes it possible to achieve the condition of sparsity and continuity of the estimation \cite{fan2001variable}.
The results in later sections can be applied to SCAD, MCP, capped-$\ell_1$, and many others.

\subsection{Directional stationary points}
\label{sec:d-stationary}

Another important definition in this paper is the d(irectional)-stationary point, which is used to describe the set of stationary solutions we are interested in this paper.
We give the definition of the d-stationary point in the following.
\begin{definition}
[\textit{d-stationary point}] A vector $\hat{\beta} \in \Omega$ is a d-stationary point to a function $Q(\beta)$ if the directional derivative, which is defined as
$$
Q'(\hat{\beta}; \beta - \hat{\beta}) = \lim_{\tau \rightarrow 0+}\frac{Q(\hat{\beta} + \tau(\beta - \hat{\beta})) - Q(\hat{\beta})}{\tau},
$$
satisfies $Q'(\hat{\beta}; \beta - \hat{\beta}) \geq 0$ for all $\beta \in \Omega$.
\end{definition}
We  prove later that under some proper conditions on the penalty function as well as on the design matrix (which in some general cases are about the loss functions), a set of  d-stationary solutions to the DC programming problem are $\sqrt{n}$ consistent estimators with a high probability (Theorem \ref{thm:estimation}).
Under further cnditions, it also recovers the true support in the unknown parameter with a high probability (Theorem \ref{thm:support}).

A motivation of choosing the directional stationary solutions (which are the directional stationary points in the corresponding optimization problem) rather than stationary solutions of other kinds, such as that of a critical point for DC programs, is provided in \cite{ahn2016difference}.
The authors argue that the directional stationary solutions are the sharpest
kind among stationary solutions of other kinds in the sense a directional stationary solution must be stationary according to other definitions of stationarity.
In the above sense, the d-stationary solutions possess minimizing properties that are not in general satisfied by stationary solutions of other kinds.
We refer to the original paper for a more detailed discussion.


\section{Formulation and assumptions}
\label{sec:Formulation}

In this section, we first give our detailed formulation in Section \ref{sec:formulate}.
We discuss the scale invariant properties of the formulation with some specific form of penalties in Section \ref{sec:scalefree}.
We then list the assumptions on the penalty functions in Section \ref{sec:hassume}, on the d-stationary solutions in Section \ref{sec:dsolassume}, on the  design matrix for our analytical study and corresponding justifications in Section \ref{sec:Xassume}.

\subsection{Formulation}
\label{sec:formulate}

We present our problem formulation in the following.
Recall that the widely-known SCAD \citep{fan2001variable} and MCP \citep{zhang2010nearly} choose their regularization term (i.e., the penalty function) as a function in the form, $\lambda|t| - h_{\lambda}(t)$, where the second term $h_{\lambda}(t)$ has a continuous first order derivative.
Motivated by the above, we propose to analyze the following parameter estimation approach:
\begin{equation}
\label{eq:obj}
\min_{\beta \in \mathbb{R}^p} \frac{1}{2n} \|Y - X\beta\|_2^2 + \lambda\|\beta\|_1 - h_{\lambda}(\beta),
\end{equation}
where function $h_{\lambda}(\beta)$ is assumed to be convex and the model is specified in (\ref{eq:LR}).
As we have shown in Table \ref{table:penalty}, popular non-convex penalties in the literature can all be expressed in DC form.
In our formulation, following the approaches in the main stream methodology, we focus on separable penalty, that is we have
\[
P(\beta) = \lambda\|\beta\|_1 - h_{\lambda}(\beta) = \sum^p_{i=1} \lambda |\beta_i| - h_{\lambda}(\beta_i),
\]
for $\beta = (\beta_1,\ldots,\beta_p) \in \mathbb{R}^p$.
Note that based on the context, function $h_{\lambda}(\cdot)$ can take both univariate and multivariate inputs.
In our formulation, the univariate function $h_{\lambda}(\cdot)$ is assumed to be convex.
Its properties are further specified  later.

\subsection{Scale-invariant property}
\label{sec:scalefree}
In real world of processing data, programmers always perform rescaling on the raw data set.
We can make our formulation scale-invariant by assuming the following format of the penalty function.
\begin{assumption}
\label{ass:scalefree}
$p_\lambda(t) = \lambda^2p(\frac{t}{\lambda})$.
\end{assumption}
Suppose we scale the model in (\ref{eq:LR}) by a scalar $c$, which leads to the following model:
\begin{equation*}
cY = X(c\beta^{\ast}) + c\epsilon,
\end{equation*}
Let $F(\beta,\lambda) = \frac{1}{2n} \|Y - X\beta\|_2^2 + \lambda\|\beta\|_1 - h_{\lambda}(\beta)$ denote the objective function, corresponding to (\ref{eq:LR}). Let $F(c\beta,c\lambda) = \frac{1}{2n} \|cY - X(c\beta)\|_2^2 + c\lambda\|c\beta\|_1 - h_{c\lambda}(c\beta)$ denote the objective function, corresponding to the scaled model.
One can easily verify that for any given positive scalar $c$,
\begin{equation}
\min_{\beta \in \mathbb{R}^p} F(c\beta,c\lambda) ,
\end{equation}
is equivalent to the original problem of $\min_{\beta \in \mathbb{R}^p} F(\beta,\lambda)$ in (\ref{eq:obj}) with scale free penalties such as SCAD, $\ell_1$, MCP, capped-$\ell_1$, which have the common form stated in Assumption \ref{ass:scalefree}. One can verify that most of the functions in Table \ref{table:penalty} satisfy this scale free condition except the logarithmic function and the transformed $\ell_1$.
However, according to the unification DCA in Section \ref{SEC:CCCP}, where in each iteration, by using only the linear approximation, we solve a re-weighted LASSO problem, which is scale invariant.

\subsection{Assumptions on $h_{\lambda}(\cdot)$}
\label{sec:hassume}

We present the assumptions that we  need in the analytical study in this section.
Recall that our penalty function has the form $P(\beta) = \lambda\|\beta\|_1 - h_{\lambda}(\beta)$.
Notice that the first term of the DC penalty is always the $\ell_1$ function.
We  specify our assumptions on the univariate function $h_{\lambda}(\beta)$, for $\beta\in\mathbb{R}$.
We  also require the regularity conditions on the design matrix $X$, which is articulated later.

The following assumptions are utilized in our analysis. We  briefly discuss the assumptions and  argue that our assumptions are equivalent or weaker to conditions in most existing work.

\begin{assumption}
\label{hderivative} We have
$\sup_{t \in \mathbb{R}}|h_{\lambda}'(t)| \leq \lambda$.
\end{assumption}
\begin{assumption}
\label{hsymmetry} We have
$h_{\lambda}(t)$ is symmetric about 0.
\end{assumption}
Both Assumption \ref{hderivative} on the non-negativity of the penalty and Assumption \ref{hsymmetry} on the symmetry of the penalty function are standard assumptions in the literature.
Assumption \ref{hderivative} makes sure that the penalty function $P_{\lambda}(\beta_i)$ is nonnegative.
In fact, we can even relax this condition to $\sup_{t \in \mathbb{R}}|h_{\lambda}'(t)| \leq \lambda$ as long as the first order derivative of the function $h_{\lambda}(t)$ is uniformly bounded in the real line.
\begin{assumption}
\label{ass:hcvxity}
$h_{\lambda}'(t)$ is monotonically increasing and there exist two nonnegative constants $\eta^- \geq \eta^+ \geq 0$ such that for any $t_2 > t_1$:
\begin{equation}
0 \leq \eta^+ \leq \frac{h_{\lambda}'(t_2) - h_{\lambda}'(t_1)}{t_2 - t_1} \leq \eta^-
\end{equation}
\end{assumption}
Regarding Assumption \ref{ass:hcvxity}, the lower bound $\eta^+$ of the convexity of the function $h_{\lambda}(t)$ is usually assumed to be $0$ in other works, such as in the SCAD and the MCP.
The upper-bound of the convexity $\eta^-$ is used to control the convexity of the function $h_{\lambda}(t)$.
If $h_{\lambda}(t)$ has a ``lot'' of convexity, we are not able to have the Restricted Strong Convexity of the objective function later.
On the other hand, this can be regarded as requiring the first order derivative of $h_{\lambda}(t)$ to be continuous.
The continuity assumption together with Assumption \ref{hderivative} and \ref{ass:zero}  ensure that Assumption \ref{ass:hcvxity} holds.

\begin{assumption}
\label{ass:horigin} We have
$h_{\lambda}(0) = h_{\lambda}'(0) = 0$.
\end{assumption}
Assumption \ref{ass:horigin} is utilized to ensure the soft thresholding property of the penalty function \cite{fan2001variable}, recalling that the singularity of the whole penalty function at $0$.

\begin{assumption}
\label{ass:zero}
For some positive $\zeta$, we have $h'_{\lambda}(t) = \lambda$ for all $|t| \geq \zeta$.
\end{assumption}

Assumption \ref{ass:zero} is based on the fact \cite{fan2001variable} that making sure $h_{\lambda}'(t) = \lambda$ for $t$ positive enough and $h_{\lambda}'(t) = -\lambda$ for $t$ negative enough  help producing an unbiased estimator.
Recall that one of the main reasons of considering a generalized version of the LASSO method is the bias of the estimation from LASSO.

Below, we make a table of the penalties discussed in Table \ref{table:penalty} and presents the decomposition to $\lambda|t| - h_{\lambda}(t)$. We also listed the properties that each of the $h_{\lambda}(t)$ holds.

\begin{table}[htp]
\begin{center}
\begin{tabular}{|c| c| c| c| c|}
\hline
Penalty & $h(t)$ & sgn(t)$h'(t)$ & Convexity  & Assumptions\\
        &              &            & measure & \\
\hline
$\ell_1$&0 &0 &0  & All except  \ref{ass:zero} \\
\hline
Capped-$\ell_1$&$\max\{0,\frac{2t}{\gamma}-1\}$ &  $\frac{2}{\gamma}I\{|t|>\gamma/2\}$& $\infty$ & All except  \ref{ass:hcvxity}\\
\hline
MCP& $(|t|- \gamma/2)I\{|t|>\gamma\}$ & $\min\{\frac{|t|}{\gamma},1\}$ & $\gamma^{-1}$ & All \\
   & $+\frac{t^2}{2\gamma}I\{|t|<\gamma\}$ & & & \\
   \hline
SCAD &$\frac{(|t|-1)^2}{2(\gamma-1)}I\{1<|t|<\gamma\} $ & $\frac{|t| - 1}{\gamma - 1} I\{1<|t|<\gamma\}$ & $(\gamma - 1)^{-1}$ &All \\
 & $+ (|t| - \frac{\gamma+1}{2})I\{|t|\geq\gamma\}$ & $+ I\{|t|\geq \gamma\} $ & & \\
\hline
Transformed $\ell_1$&$|t| - \frac{a|t|}{a+ |t|}$ & $\frac{(a+ |t|)^2 - a^2}{(a+ |t|)^2} $&$\frac{2}{a}$ &All except  \ref{ass:zero} \\
\hline
Logarithmic&$|t| - \log(|t|+1)$ &  $\frac{|t|}{|t|+1}$& $1$& All except  \ref{ass:zero}\\
\hline
\end{tabular}
\end{center}
\caption{The penalties in the sparse estimation literature and their properties with respect to our assumptions.
The first column gives the name of the methods.
The second column presents the $h$-function, which is the second component in the DC decomposition ($p=g-h$) of the corresponding penalty function.
The third column contains their first derivatives on the positive axe.
This is to verify Assumption \ref{hderivative}.
The fourth column computes for the quantities that are raised in Assumption \ref{ass:hcvxity}.
The last column summarizes the assumptions that the corresponding penalty satisfies. }
\label{table:decomp}
\end{table}

From  Table \ref{table:decomp}, we can see that, SCAD and MCP penalty class satisfy all of the assumptions.
While Capped-$\ell_1$ has discontinuous first order derivative, which violates the Assumption \ref{ass:hcvxity}.
In order to extend the theories in this work to Capped-$\ell_1$ penalty, performing smoothing around the non-differentiable point is enough.
We re-scaled the linear term in Transformed $\ell_1$ to match the assumptions.
$\epsilon$ in Logarithmic penalty is chosen to be 1 in the Table \ref{table:decomp}.

\subsection{Assumptions on the d-stationary solution}
\label{sec:dsolassume}

Besides the assumptions on the penalty functions, we also list the assumptions necessary for the loss function. These are about the design matrix in case of linear model with the least square loss function.
\begin{assumption}
\label{nonnegativity}
Let $\beta^{\ast}$ be the unknown true parameter, $\hat{\beta}$ be the d-stationary solution to problem (\ref{eq:obj}), which satisfies the following condition:
\[
\frac{1}{n}X_j^TX(\beta  - \hat{\beta})\mbox{sign}(\hat{\beta}_j) \geq c\lambda,\text{ for all $j \in S^c$, $\beta = \beta^{\ast}$ with $c \in (0,1)$.}
\]
\end{assumption}
\begin{remark}
The above Assumption \ref{nonnegativity} is no stronger than the assumptions used in LASSO estimator \cite{wainwright2009sharp}.
We show below that in the proof of LASSO estimator, it corresponds to when $c = \frac{1}{2}$.
Let $\hat{\beta}^{lasso} = \arg \min_{\beta \in \mathbb{R}^p} \frac{1}{2n} \|Y - X\beta\|_2^2 + \lambda\|\beta\|_1$. Recall that in \cite{wainwright2009sharp}, by the First Order Condition (FOC) at $\hat{\beta}^{lasso}$, we have
\[-\frac{1}{n}X^T(Y - X\hat{\beta}^{lasso}) + \lambda \partial \|\hat{\beta}^{lasso}\|_1 = 0,\]
where $\partial \|\hat{\beta}^{lasso}\|_1$ is a subgradient at $\hat{\beta}^{lasso}$ for $\|\beta\|_1$. Multiply by $(\beta^{\ast} -\hat{\beta}^{lasso})^T$ on both side, we have
\[-\frac{1}{n}(\beta^{\ast} -\hat{\beta}^{lasso})^TX^T(Y - X\hat{\beta}^{lasso}) + \lambda (\beta^{\ast} -\hat{\beta}^{lasso})^T\partial \|\hat{\beta}^{lasso}\|_1 = 0.\]
Since
\begin{equation}
\begin{split}
&(\beta^{\ast} -\hat{\beta}^{lasso})^T\partial \|\hat{\beta}^{lasso}\|_1\\
=& (\beta^{\ast} -\hat{\beta}^{lasso})_S^T\partial \|\hat{\beta}^{lasso}_S\|_1
+ (\beta^{\ast} -\hat{\beta}^{lasso})_{S^c}^T\partial \|\hat{\beta}^{lasso}_{S^c}\|_1 \\
=& (\beta^{\ast} -\hat{\beta}^{lasso})_S^T\partial \|\hat{\beta}^{lasso}_S\|_1
- \|\hat{\beta}^{lasso}_{S^c}\|_1
\end{split}
\end{equation}
Plugging into the FOC, we have
\begin{equation}
\begin{split}
&\frac{1}{n}(\beta^{\ast} -\hat{\beta}^{lasso})^TX^T(Y - X\hat{\beta}^{lasso}) \\
=& \frac{1}{n}(\beta^{\ast} -\hat{\beta}^{lasso})^TX^TX(\beta^{\ast} -\hat{\beta}^{lasso}) +  \frac{1}{n}(\beta^{\ast} -\hat{\beta}^{lasso})^TX^T\epsilon\\
=& \frac{1}{n}(\beta^{\ast} -\hat{\beta}^{lasso})^TX^TX(\beta^{\ast} -\hat{\beta}^{lasso}) + \frac{1}{n}(\beta^{\ast} -\hat{\beta}^{lasso})_S^TX_S^T\epsilon\\
 &+ \frac{1}{n}(\beta^{\ast} -\hat{\beta}^{lasso})_{S^c}^TX_{S^c}^T\epsilon\\
 =& \frac{1}{n}(\beta^{\ast} -\hat{\beta}^{lasso})^TX^TX(\beta^{\ast} -\hat{\beta}^{lasso}) + \frac{1}{n}(\beta^{\ast} -\hat{\beta}^{lasso})_S^TX_S^T\epsilon\\
 &- \sum_{S^c}\frac{1}{n}|\hat{\beta}^{lasso}|_iX_i^T\epsilon\text{sign}(\hat{\beta}^{lasso}_i)\\
=& \lambda(\beta^{\ast} -\hat{\beta}^{lasso})_S^T\partial \|\hat{\beta}^{lasso}_S\|_1
- \lambda\|\hat{\beta}^{lasso}_{S^c}\|_1
\end{split}
\end{equation}
where in the first equality, we plugged in $Y = X\beta^{\ast} + \epsilon$. If we have the condition that $\frac{1}{n}X_i^TX(\beta^{\ast} -\hat{\beta}^{lasso})\text{sign}(\hat{\beta}^{lasso}_i) > c\lambda$ ($c = \frac{1}{2}$ in LASSO) for all $i \notin S$, we  have
\begin{equation}
\begin{split}
&\frac{1}{n}(\beta^{\ast} -\hat{\beta}^{lasso})^TX^TX(\beta^{\ast} -\hat{\beta}^{lasso}) \\
&\leq \frac{3}{2}\lambda\|\hat{\beta}^{lasso}_{S}\|_1
- c\lambda\|\hat{\beta}^{lasso}_{S^c}\|_1
\end{split}
\end{equation}
with high probability. Similarly, we made Assumption \ref{nonnegativity} in the generalized penalized regression to ensure good property of the solution.
\end{remark}

\begin{remark}
Since the condition in Assumption \ref{nonnegativity} cannot be verified directly, in real data analysis, we can use the following checkable conditions instead:
\[
\frac{1}{n}X_j^T(Y  - X\hat{\beta})\mbox{sign}(\hat{\beta}_j) \geq c\lambda,\text{ for all $j \in S^c$ such that $\hat{\beta}_j \neq 0$, $\beta = \beta^{\ast}$,}
\]
where $c$ is defined in Assumption \ref{nonnegativity}.
If the above holds, the Assumption \ref{nonnegativity}  hold with high probability using similar argument of sub-Gaussian random variables as in the proof of Theorem \ref{cor:Bon}.
\end{remark}

\subsection{Assumptions on the design matrix $X$}
\label{sec:Xassume}

\begin{definition}
The restricted strong convexity (RSC) condition on model matrix $X$ with respect to $\mathcal{C}$ is the following, there exists some constant $\gamma > 0$ such that:
\[\frac{\frac{1}{N}\nu X^TX \nu}{\|\nu\|_2^2} \geq \gamma\text{ for all nonzero } \nu \in \mathcal{C}\]
\end{definition}
\noindent here $\gamma$ is called the restricted eigenvalue bound with regard to $\mathcal{C}$.
\begin{assumption}
\label{RSC}
Denote $C_{\mathcal{S}}$ by the diagonal matrix with $\{c_i, i \in \mathcal{S}\}$, $C_{\mathcal{S}^C}$ by the diagonal matrix with $\{c_i, i \in \mathcal{S}^C\}$, the restricted eigenvalues (RE) condition holds on the following set with some positive $c$ defined in Assumption \ref{nonnegativity}:
\[
\mathcal{C} = \left\{\nu \in \mathbb{R}^p \left| \|C_{\mathcal{S}^C} \cdot \nu_{\mathcal{S}^C}\|_1
\leq \frac{5}{2c}\|C_{\mathcal{S}} \cdot \nu_{\mathcal{S}}\|_1
\right.
\right\},
\]
where $\cdot$ indicates the matrix-vector multiplication.
\end{assumption}
\noindent We have $\mathcal{C} \subset \mathbb{R}^p$ strictly since it is of the form of cone.

The RSC (Assumption \ref{RSC}) is a standard assumption in the literature for proving the consistency results of regularized high-dimensional sparse estimation problems.

\section{Consistency results for some d-stationary solutions}\label{MAINRESULTS}
We  prove our main results in this section.
The non-asymptotic upper bound for estimation errors is derived in Section \ref{subsec:non-asymptotic}.
As a corollary, we provide the upper bound for prediction errors as a byproduct  in Section \ref{sec:predictionError}.
We provide the results regarding asymptotic consistency of the estimation in Section \ref{sec:asymptoticConsistency}.
The asymptotic consistency in support recovery is discussed in Section \ref{subsec:support-recovery}.
\subsection{Non-asymptotic upper bound for estimation errors}\label{subsec:non-asymptotic}

In this section, we  present our results on the  non-asymptotic upper bound for the estimation error.
We mainly use the assumptions on the model matrix to prove that the difference between the ground truth $\beta^{\ast}$ and the d-stationary solution  $\hat{\beta}^{\lambda}_{n}$ will be in a cone-like set, where  we have the restricted strong convexity (RSC) assumption hold (as defined in Assumption \ref{RSC}).
Without  Assumption \ref{ass:hcvxity} on the continuity of the first order derivative on the function $h_{\lambda}(t)$, we will be able to obtain the upper-bound of the $\ell_2$ distance between the ground truth and the d-stationary estimation.

\begin{theorem}
\label{thm:estimation}
Suppose $h_{\lambda}(t)$ satisfies Assumptions \ref{hderivative}, \ref{hsymmetry}, \ref{ass:horigin}, design matrix $X$ satisfies the restricted strong convexity with respect to $\mathcal{C}$, which is defined in Assumption \ref{RSC} with $c_i = 1$ for $i = 1, \cdots, p$,  with $\lambda \geq \frac{2\|X^T\epsilon\|_{\infty}}{n} $.
If we further assume Assumption \ref{nonnegativity} holds at the d-stationary solution $\hat{\beta}^{\lambda}_n$, we will have the upper bound for estimation error on $\beta^{\ast}$ with the d-stationary estimation $\hat{\beta}^{\lambda}_{n}$ as $n \rightarrow \infty$:
\[
\|\hat{\beta}^{\lambda}_n - \beta^{\ast}\|_2 \leq  \frac{5}{2\gamma}\lambda\sqrt{\|\beta^{\ast}\|_0}
 \]
\end{theorem}

The proof of the above Theorem \ref{thm:estimation} is delayed to Appendix \ref{proof:estimation}.
The results in Theorem \ref{thm:estimation} suggest that the d-stationary solution to  Problem (\ref{eq:obj}), under mild conditions, will be able to retrieve the information in the unknown parameter $\beta^{\ast}$ with error bounded within the order of $O(\frac{\lambda \sqrt{|S|}}{\gamma})$, which is optimal.

\subsection{Non-asymptotic upper bound for prediction errors}\label{sec:predictionError}
From the proof of Theorem \ref{thm:estimation}, we will be able to further give the upper bound for the prediction error below.

\begin{corollary}
\label{cor:estimation}
Under the assumptions of Theorem \ref{thm:estimation}, we can further get the upper-bound for the prediction error as:
\[
\left\|\frac{1}{\sqrt{n}}X(\beta^{\ast}  - \hat{\beta}^{\lambda}_n)\right\|^2_2
\leq (\frac{5}{2}\lambda)^2 \frac{1}{\gamma}\sqrt{|S|}\]
\end{corollary}
The proof of Corollary \ref{cor:estimation} is straight forward according to the proof in Theorem \ref{thm:estimation} and is postponed to Appendix \ref{proof:corestimation}.

\subsection{Asymptotic convergence rate}\label{sec:asymptoticConsistency}
If we further assume that the errors are independent sub-Gaussian distributed, we will be able to bound the estimation error in the order of $\sqrt{\frac{ |S|\log{p}}{n}}$ with high probability.

\begin{corollary}
\label{cor:Bon}
Under the assumptions of Theorem \ref{thm:estimation}, if we further assume that the errors are from independent sub-Gaussian with variance parameter $\sigma^2$ and mean $0$, we will have the following hold with probability at least $1 - 2\exp{(-\frac{\tau - 2}{2}\log{p})}$
\[\|\hat{\beta}^{\lambda}_n - \beta^{\ast}\|_2 \lesssim \frac{5}{\gamma}\sigma\sqrt{\frac{\tau |S|\log{p}}{n}}  .\]
\end{corollary}
We provide the proof of Corollary \ref{cor:Bon} in Appendix \ref{proof:corBon}.
\begin{remark}
All the results above are considering the problem in (\ref{eq:obj}).
However, the conclusions will still hold for the constrained version of problem (\ref{eq:obj}) as long as the assumptions are satisfied.
The results in this work assumes the d-stationary solution that satisfies our assumptions exists.
We will justify the existence of the d-stationary solution satisfying our assumptions in Section \ref{sec:existence}.
\end{remark}
\subsection{Support recovery}\label{subsec:support-recovery}

In this section, we will first provide the KKT conditions for d-stationary solutions, which says that the d-stationary condition in our work is equivalent to the first order condition in case of no constraints.
Then we  prove the Restricted Strong Convexity for the objective function in Problem \eqref{eq:obj} under some regularity conditions.
By usage of the oracle estimator defined later in Problem \eqref{estoracle}, we will be able to prove the support recovery consistency of some of the d-stationary solutions to Problem \eqref{eq:obj}.

\begin{lemma}
\label{FOClemma}
Let $F(\beta) = L(\beta) + g(\beta) - h(\beta)$, where $L(\beta)$, $g(\beta)$, $h(\beta)$ are convex with $\beta \in \mathbb{R}^p$. Further assume that $L(\beta)$ and $h(\beta)$ have continuous first order derivative, $g(\beta) = \|\beta\|_1$. Let $\beta_0$ be a d(irectional)-stationary solution to $F(\beta)$, we  have the following first order condition (FOC) hold at $\beta_0$. We will be able to get the following equivalent condition: $\beta_0$ is a d(irectional)-stationary solution to $F(\beta)$ if and only if there exists some $z \in \partial g(\beta_0)$, where  $\partial g(\beta_0)$ is the set of subgradient of $g(\beta)$ at $\beta_0$, such that:
\begin{equation}
 \nabla L(\beta_0) +  z -  \nabla h(\beta_0) = 0,
\end{equation}
 where $ \nabla L(\beta_0)$,  $\nabla h(\beta_0)$ is the gradient of $L$, $h$ at $\beta_0$.
\end{lemma}

The above Lemma \ref{FOClemma} states the equivalence between d-stationary solution and first order condition (FOC) in the unconstrained case. While in constrained case, this does not necessarily hold.
From the proof Lemma \ref{FOClemma}, we can derive similar conditions for ``local maximals'' for $\tilde{\beta}$.
We obtain that as long as $\min_{i = 1}^p\{\tilde{\beta}_{i}\} = 0$, it will only satisfy the condition for ``local'' minimals and thus be a d-stationary solution.
Thus, in order to find the d-stationary solution, we only need to find a $\beta_0$ such that, there exists a vector $z \in  \partial\|\beta_0\|_1$, the subgradient of function $\|\beta\|_1$ at $\beta = \beta_0$, such that $\nabla L(\beta_0) - \nabla h(\beta_0) + z = 0$. Furthermore, if $\min_{i = 1}^p\{z_{i}\} < 1$, which is known as the strict dual feasibility condition \citep{wainwright2009sharp}, it will be satisfying the condition for ``local'' minimals.

\begin{remark}
Generally, a d-stationary solution is not necessarily local minimal. For example, for a differentiable function $f(x,y) = x^2 - y^2$, where both the function $g(x,y)$ and $h(x,y)$ are differentiable (slightly different from the above situation in Lemma \ref{FOClemma}), at a saddle point $(0,0)$, which is stationary with $\mathbf{0}$ gradient, the directional derivative at this point will all be $0$, which makes it a d-stationary solution however not a local minimal. Another example would be $f(x,y) = |x| - y^2$ at the saddle point $(0,0)$.
\end{remark}
\begin{remark}
The necessary condition to be a local minimal is being a d-stationary point in the feasible region.
\end{remark}

The following Lemma shows the RSC of the Problem \eqref{eq:obj}.
\begin{lemma}
\label{lem:cvxity}
Under Assumption \ref{RSC} with $h_\lambda$ satisfying Assumptions \ref{hderivative}, \ref{hsymmetry}, \ref{ass:hcvxity}, \ref{ass:horigin}, let $\beta_1$, $\beta_2 \in \mathbb{R}^p$ such that $\nu = \beta_1 - \beta_2 \in \mathcal{C}$, where $\mathcal{C}$ is defined in Assumption \ref{RSC}. Then $f_{\lambda}(\beta) = \frac{1}{2n}\|Y - X\beta\|^2_2 - h_{\lambda}(\beta)$ will satisfy the restricted strong convexity given that $\gamma > \eta^-$:
\begin{equation}
\label{cvxity}
f_{\lambda}(\beta_2) \geq f_{\lambda}(\beta_1) + \nabla f_{\lambda}(\beta_1)^T(\beta_2 - \beta_1) + \frac{\gamma - \eta^-}{2}\|\beta_2 - \beta_1\|^2_2
\end{equation}
\end{lemma}
The proof can be found in Appendix \ref{proof:lemcvxity}
\subsubsection{Oracle estimator}
The oracle estimator is defined as follows:
\begin{equation}
\label{estoracle}
\beta^O = \arg\min_{\beta \in \mathbb{R}^p, \beta_{S^c} = 0} \frac{1}{2n}\|Y - X\beta\|^2_2.
\end{equation}
The oracle estimator is obtained as if there is an oracle telling the true support of the underlying unknown estimator. According to the definition of oracle estimator $\beta^O$, we are able to provide the following $\ell_{\infty}$ error bound between $\beta^O$ and $\beta^{\ast}$.
We also demonstrate that $\beta^O$ is a d-stationary solution to the DC-penalized Problem (\ref{eq:obj}), which we are interested in this paper. The following Theorem \ref{oracle} and Lemma \ref{oracledstat} also appeared in the work by Wang et al. \cite{wang2014optimal}

\begin{theorem}
\label{oracle}
Under Assumption \ref{RSC}, the oracle estimator is the unique global minimizer of (\ref{estoracle}).
If the noise is independent sub-Gaussian with variance parameter $\sigma^2$, the oracle estimator will satisfy the following $\ell_{\infty}$ error bound with high probability.
\[\|\beta^O - \beta^{\ast}\|_{\infty} \leq C\sigma\sqrt{2/\gamma}\sqrt{\frac{\log s}{n}}.\]
\end{theorem}
The proof is in Appendix \ref{proof:oracle}.

\begin{lemma}
\label{oracledstat}
Under Assumption \ref{RSC} with $h_\lambda$ satisfying Assumptions \ref{hderivative}, \ref{hsymmetry}, \ref{ass:hcvxity}, \ref{ass:horigin}, \ref{ass:zero}, let $\beta^O$ be the aforementioned oracle estimator.
Assume further that for the ground truth $\beta^\ast$, we have $\min_{i = 1}^p|\beta_i^{\ast}| > 2\zeta$, for  $\zeta > 0$.
There exists a subgradient $\xi^O \in \partial \|\beta^O\|_1$ (where $\partial \|\beta^O\|_1$ stands for the subgradient of function $\|\beta\|_1$ at $\beta = \beta^O$) such that for any $\beta \in \mathbb{R}^p$:
\begin{equation}
\label{focOracle}
 (\beta - \beta^O)^T \left(\nabla f_{\lambda}(\beta^O) + \lambda \xi^O\right) \geq 0
\end{equation}
\end{lemma}

The above Lemma \ref{oracledstat} assumes that the penalty on the parameters will be a constant when the parameters are  large. As it requires Assumption \ref{ass:zero}, the result is not applicable for transformed $\ell_1$ and logarithmic penalties. We postpone the proof to Appendix \ref{proof:oracledstat}.

\begin{lemma}
\label{subset}
Under the assumptions in Lemma \ref{oracledstat}, let $\hat{\beta}$ be a d-stationary solution to (\ref{eq:obj}) satisfying Assumption \ref{nonnegativity}, and $\beta^O$ be the oracle estimator. The following will hold with large probability:
\begin{equation}
\nu = \hat{\beta} - \beta^O \in \mathcal{C}.
\end{equation}
\end{lemma}

The proof is in Appendix \ref{proof:subset}.
Based on the previous results of the oracle estimator $\beta ^O$ and properties of the d-stationary estimator $\hat{\beta}$, we will now be able to prove the support recovery consistency for our generalized DC-penalized model.

\begin{theorem}
\label{thm:support}
Under the conditions of Lemma \ref{subset}, we will have $supp(\hat{\beta}) = supp(\beta^O) = supp(\beta^{\ast})$ with high probability.
\end{theorem}
The proof is provided in Appendix \ref{proof:thmsupport}.
We prove the support recovery consistency for a set of d-stationary solutions, it implies that a set of the convergence points (satisfying Assumption \ref{nonnegativity}) from the DCA will converge to the oracle estimator which is unique.
In the work from Wang et al, they prove that the convergence point from each stage of the specific algorithm converges to the oracle estimator in the linear model setting.
The above results also inform us how we should choose the penalty function such that the d-stationary solution will be support recovery consistent. The penalty needs to be a constant when the parameter gets larger (Assumption \ref{ass:zero}), so that the resulting oracle estimator will be a d-stationary solution to the original Problem (\ref{eq:obj}). Assumption \ref{ass:hcvxity} is necessary for the restricted strong convexity in $\mathcal{C}$.

\section{DC penalty with generalized loss functions}\label{GLOSS}

In the previous section, we mainly focus on the linear model scenario.
Most of the analysis can be readily extended to generalized loss functions such as logistic loss function, etc.
Below, we will present the formulation of DC penalized likelihood and provide the statistical analysis regarding to the d-stationary solutions.

We begin with a brief review on the exponential family. Exponential family is a family of distributions with the  probability density proportional to $P(Y|X,\beta^\ast) \propto \exp\{\frac{YX^T\beta^\ast - \psi(X^T\beta^\ast)}{c(\sigma)}\}$, where $c(\sigma)$ is a scaling parameter and $\psi(\cdot)$ is the cumulant function.
According to \cite{lehmann1998theory}, one standard property of exponential family is
$$\psi'(X^T\beta^\ast) = \mathbb{E}[Y|X,\beta^\ast,\sigma].$$
Given that $\psi(\cdot)$ is a univariate convex function, let $L(\beta) = \psi(X^T\beta) - YX^T\beta$ be the negative log likelihood function, $L_n(\beta) = \frac{1}{n} \sum_{i = 1}^n (\psi(X_i^T\beta) - Y_iX_i^T\beta)$ be the sample average of the negative log likelihood function, one can easily check that $\mathbb{E}[\nabla L(\beta^\ast)] = 0$ and $\nabla^2 L_n(\beta) \geq 0$.
This implies that $L_n(\beta)$ is convex.
In the following, we might omit the subscript $n$ in the expression of $L_n(\beta)$ where no confusion will rise.
In order to estimate the sparse ground truth $\beta^\ast$, we will solve the following DC penalized optimization problem:
\begin{equation}
\label{eq:objGLM}
\min_{\beta \in \mathbb{R}^p} \frac{1}{n} \sum_{i = 1}^n (\psi(X_i^T\beta) - Y_iX_i^T\beta) + \lambda\|\beta\|_1 - h_{\lambda}(\beta),
\end{equation}

Below, we will state the assumptions on the generalized loss functions, which enable us to provide the analysis that the error between the  d-stationary solution $\hat{\beta}$ and the ground truth $\beta^\ast $ is of the order $\mathbb{O}(\frac{17\lambda}{8(\gamma - \eta^-)}\sqrt{|S|})$.
\begin{assumption}
\label{ass:gradientLoss}
Let $\beta^\ast$ represent the ground truth of the unknown parameter, $L(\beta)$ be the negative log likelihood function. Assume that the infinity norm of the gradient of the loss function at the ground truth $\|\nabla L(\beta^\ast)\|_\infty \leq \frac{\lambda}{8}$.
\end{assumption}

\begin{assumption}
\label{ass:data}
Let $\beta^{\ast}$ be the unknown true parameter, $\hat{\beta}$ be the d-stationary solution to problem (\ref{eq:obj}), which satisfies the following condition:
\[
\|\nabla h_\lambda(\hat{\beta}_{S^c})\|_{\infty} \leq (1-c)\lambda,\text{ with $c \in (0,1)$.}
\]
\end{assumption}

\begin{assumption}
\label{ass:RSCGLM}
Let $\beta^\ast$ represent the ground truth of the unknown parameter,  $L(\beta)$ be the negative log-likelihood function.
Assume that the following restricted strong convexity holds on the set $\mathcal{C}$,
\begin{eqnarray*}
\mathcal{C} &=& \left\{\nu \in \mathbb{R}^p \left| \|C_{\mathcal{S}^C}\nu_{\mathcal{S}^C}\|_1
\leq \frac{4+c}{c}\|C_{\mathcal{S}}\nu_{\mathcal{S}}\|_1
\right.
\right\},  \\
L(\beta_1) &\geq & L(\beta_2) + \nabla L(\beta_2)^T(\beta_1 - \beta_2) + \frac{\gamma}{2}\|\beta_1 - \beta_2\|^2_2,
\end{eqnarray*}
for any $\beta_1$ and $\beta_2$ such that $\beta_1-\beta_2 \in \mathcal{C}$.
\end{assumption}

\begin{theorem}
\label{thm:gloss}
Let $\hat{\beta}$ be the d-stationary solution to the penalized loss function in (\ref{eq:objGLM}).
Suppose $h_{\lambda}(t)$ satisfies Assumptions \ref{hderivative}, \ref{hsymmetry}, \ref{ass:hcvxity}, \ref{ass:horigin},  assume further that Assumptions \ref{ass:gradientLoss}, \ref{ass:RSCGLM} and \ref{ass:data} hold,
if $\gamma > \eta^-$,
$$
\left\|\frac{1}{n} \sum_{i = 1}^n(\psi'(X_i^T\beta^\ast)X_i - Y_iX_i)\right\|_{\infty} \leq \frac{c}{2},
$$
we will have the following upper-bound for the estimation error of the d-stationary solution
$$\|\beta^\ast - \hat{\beta}\|_2 \leq \frac{17\lambda}{8(\gamma - \eta^-)}\sqrt{|S|}$$
\end{theorem}
The proof is provided in Appendix \ref{proof:thmgloss}.

\subsection{Existence of d-stationary solution}\label{sec:existence}
In this section, we will show the existence of the d-stationary solutions we studied above. It is easy to see that in the linear regression setting with square loss, the oracle estimator is a d-stationary solution under  suitable conditions we stated in Lemma \ref{oracledstat}. For general settings with generalized loss functions, let $r_0 > 0$ be such that $h_\lambda'(r_0) = (1-c)\lambda$, consider the following constrained problem:
\begin{equation}
\label{eq:obj2}
\min_{\|\beta - \beta^\ast\|_2 \leq r} L_n(\beta) + \lambda\|\beta\|_1 - h_{\lambda}(\beta),
\end{equation}
where $r = c\lambda\sqrt{|S|}\wedge r_0$. It is straightforward to check that the sationary solutions to Problem (\ref{eq:obj2}) satisfies all the assumptions of the d-stationary solution  studied in Section \ref{MAINRESULTS} and Section \ref{GLOSS}, which verifies the existence of the wanted d-stationary solutions.

\section{Numerical Approach to Find the d-stationary Points}\label{SEC:CCCP}

In this section, we will review the efficient algorithms in the DC-literature, for finding the  local optima in the statistics and optimization areas.
 This provides a comprehensive summary on solving DC programming.
Up to our knowledge, the most classic algorithm is the Difference-of-Convex Algorithm (DCA)  studied in \cite{tao1997convex, horst1999dc, tao1997convex, sriperumbudur2012proof}, which iterates between the primal problem and the dual problem to find the local minima.
Given the DC problem below:
\begin{equation}
\label{eq:DCgeneral}
\min_{x \in \mathbb{R}^n} f(x) = g(x) - h(x),
\end{equation}
where $g(\cdot)$ and $h(\cdot)$ are convex functions.
For a function $g: \mathbb{R}^n \rightarrow \mathbb{R}$, let $g^\ast(y)$ be its convex conjugate function, which is defined as $g^\ast(y) = \sup\{x^Ty - g(x): x \in \mathbb{R}^n\}$.
We have
\begin{equation}
\begin{split}
&\inf_{x \in \mathbb{R}^n} {f(x) = g(x) - h(x)}\\
=&\inf_{x}\{g(x) - \sup_y\{x^Ty - h^\ast(y)\}\}\\
=&\inf_{x}\{ \inf_y\{g(x) +h^\ast(y) - x^Ty\}\}\\
=&\inf_{y}\{ -\sup_x\{ x^Ty- g(x) \}+h^\ast(y)\}\\
=&\inf_{y}\{h^\ast(y) - g^\ast(y)\}.\\
\end{split}
\end{equation}
Thus, by iterating between the primal and the dual problems, the  DCA will converge to a d-stationary solution.
Below shows the DCA.

\begin{algorithm}
\label{alg:DCA}
\caption{Difference-of-Convex Algorithm (DCA)}
\begin{algorithmic}[1]
\STATE Choose the initial $x_0$
\STATE \emph{loop}:
\FOR {$k \in \mathbb{N}$}
\STATE Choose $y_k \in \partial h(x_k)$.
\STATE Choose $x_{k+1} \in \partial g^\ast(y_k)$.
\IF{$(\min\{|(x_{k+1} - x_{k})_i|,|\frac{(x_{k+1} - x_{k})_i}{(x_{k})_i}|\} \leq \delta)$}
\RETURN $x_{k+1}$
\ENDIF
\ENDFOR
\end{algorithmic}
\end{algorithm}
According to \cite{tao2005dc} in Section 2.5, DCA has linear convergence rate for general DC programmings.
While in the statistics literature, Local Linear Approximation (LLA) in \cite{zou2008one} is widely used for solving regularized estimation problems with non-convex penalties.
The update at each iteration takes the LLA of the penalty function:
$$x^{k+1} = \arg\min \{g(x) - \partial h(x_k)^Tx\},$$
which is exactly the same procedure as shown in the Algorithm \ref{alg:DCA}.

In the setting of this paper, the objective is defined in (\ref{eq:obj}), where we are minimizing the objective function over all $\beta \in \mathbb{R}^p$ with the first part of the DC function as $g(\beta) = L_n(\beta) + \lambda\|\beta\|_1$, and the second part of DC function $h(\beta) = h_{\lambda}(\beta)$.
The DCA can be simplified to Local Linear Approximation (LLA) in the general case as in \cite{zou2008one}, the detailed procedures can be found in \cite{sriperumbudur2012proof}.
Specifically, if  $h(x)$ is differentiable,  we will have the following equivalent algorithm as DCA:

\begin{algorithm}
\caption{DCA (LLA)}\label{alg:DCALLA}
\begin{algorithmic}[1]
\STATE Choose the initial $\beta_0$
\STATE \emph{loop}:
\FOR {$k \in \mathbb{N}$}
\STATE Choose $z_k \in \nabla h(\beta_k)$.
\STATE $\beta_{k+1} = \arg \min L_n(\beta) + \lambda\|\beta\|_1 - \langle \beta, \nabla h(\beta_k) \rangle$.
\IF{$(\min\{|(\beta_{k+1} - \beta_{k})_i|,|\frac{(\beta_{k+1} - \beta_{k})_i}{(\beta_{k})_i}|\} \leq \delta)$}
\RETURN $\beta_{k+1}$
\ENDIF
\ENDFOR
\end{algorithmic}
\end{algorithm}

According to \cite{yuille2003concave}, DCA is exactly the formulation of Convex Concave Procedure (CCCP), which is also discussed in \cite{pang2016computing}.
Thus, under proper conditions, all results in \cite{yuille2003concave} can be applied to the problem studied here.
Since our formulation \eqref{eq:obj} is a special form of the model considered in \cite{pang2016computing}, which adopts the classical algorithm DCA (Difference-of-Convex Algorithm) in \cite{tao1997convex, horst1999dc, tao1997convex, sriperumbudur2012proof} and solves a strictly convex problem at each iteration, it is guaranteed to converge quickly to a d-stationary solution.
Since the penalty is a function of the absolute value of the estimator, one minor change to the above algorithm would be solving the following transformed optimization problem within each iteration:
\begin{equation}
\label{eq:LLA}
\beta_{k+1} = \arg \min L_n(\beta) + \sum_{i = 1}^p (\lambda - h'(|\beta_{ki}|)) |\beta_i|,
\end{equation}
which is exactly the formulation of weighted LASSO estimator and can be solved efficient using the LARS algorithm in \cite{efron2004least}.
\begin{lemma}
\label{lemma:Decrease}
By updating the parameter $\beta$ as in Procedure \ref{eq:LLA}, the objective function $F(\beta)$ defined in \ref{eq:obj} is monotonically decreasing.
\end{lemma}

In the one-step LLA procedure \cite{zou2008one}, the authors prove that starting from the maximum likelihood estimator (MLE), after one step of the LLA update, the resulting estimator is consistent when SCAD penalty function is used.
While in \cite{fan2014strong}, they prove that from the LASSO initialization, with high probability that the LLA converges to the oracle estimator in 2 iterations.
The above results can be similarly extended to our DC setting.

\section{Conclusions}\label{sec:conclusion}
In this work, we close the gap between the statistics and optimization by finding a set of d-stationary solutions to the DC penalized loss functions.
Specifically, we relax the assumptions used in \cite{ahn2016difference} and provide stronger statistical results on the penalized estimation problem.
We prove that a certain subset of d-stationary solutions in an optimization problem (with a DC objective) has the ideal statistical properties: asymptotic estimation consistency, asymptotic model selection consistency, asymptotic efficiency under linear model and the GLM settings.
We also provide the non-asymptotic upper bound for the estimation errors in both scenarios.
We unify the framework of non-convex penalized high-dimensional sparse estimation problems and the existing popular algorithms to solve the problems in a DC framework.

Several open questions remain, which might be interesting directions for future research. Since in this work, we mainly consider the unconstrained DC programming, it is unclear whether a proper constraint, which might depend on specific problems, will ensure a better set of solution or possibly unique solution to the high-dimensional sparse estimation problem.
Another direction would be more general loss functions. When the observations have outliers or missing values, it would be desiring to obtain theoretical guarantees on the sparse estimations with possibly non-convex loss functions, such as Huber loss, Cauchy loss, etc.

\bibliographystyle{imsart-number}
\bibliography{DC_LASSO}

\begin{thebibliography}{41}

\bibitem{ahn2016difference}
\begin{barticle}[author]
\bauthor{\bsnm{Ahn},~\bfnm{Miju}\binits{M.}},
  \bauthor{\bsnm{Pang},~\bfnm{Jong-Shi}\binits{J.-S.}} \AND
  \bauthor{\bsnm{Xin},~\bfnm{J}\binits{J.}}
(\byear{2016}).
\btitle{Difference-of-convex learning I: Directional stationarity, optimality,
  and sparsity}.
\bjournal{SIAM Journal on Optimization, revision under review (as of February
  2017)}.
\end{barticle}
\endbibitem

\bibitem{aleksandrov1950surfaces}
\begin{binproceedings}[author]
\bauthor{\bsnm{Aleksandrov},~\bfnm{AD}\binits{A.}}
(\byear{1950}).
\btitle{Surfaces represented as a difference of two convex functions, Russian
  Acad. Sci}.
In \bbooktitle{Dokl. Math}
\bvolume{1}.
\end{binproceedings}
\endbibitem

\bibitem{bickel2009simultaneous}
\begin{barticle}[author]
\bauthor{\bsnm{Bickel},~\bfnm{Peter~J}\binits{P.~J.}},
  \bauthor{\bsnm{Ritov},~\bfnm{Ya'acov}\binits{Y.}} \AND
  \bauthor{\bsnm{Tsybakov},~\bfnm{Alexandre~B}\binits{A.~B.}}
(\byear{2009}).
\btitle{Simultaneous analysis of {Lasso} and {Dantzig} selector}.
\bjournal{The Annals of Statistics}
\bpages{1705--1732}.
\end{barticle}
\endbibitem

\bibitem{buhlmann2011statistics}
\begin{bbook}[author]
\bauthor{\bsnm{B{\"u}hlmann},~\bfnm{Peter}\binits{P.}} \AND \bauthor{\bsnm{Van
  De~Geer},~\bfnm{Sara}\binits{S.}}
(\byear{2011}).
\btitle{Statistics for high-dimensional data: methods, theory and
  applications}.
\bpublisher{Springer Science \& Business Media}.
\end{bbook}
\endbibitem

\bibitem{candes2007dantzig}
\begin{barticle}[author]
\bauthor{\bsnm{Candes},~\bfnm{Emmanuel}\binits{E.}} \AND
  \bauthor{\bsnm{Tao},~\bfnm{Terence}\binits{T.}}
(\byear{2007}).
\btitle{The {Dantzig} selector: Statistical estimation when $p$ is much larger
  than $n$}.
\bjournal{The Annals of Statistics}
\bpages{2313--2351}.
\end{barticle}
\endbibitem

\bibitem{chen1995examples}
\begin{binproceedings}[author]
\bauthor{\bsnm{Chen},~\bfnm{Scott}\binits{S.}} \AND
  \bauthor{\bsnm{Donoho},~\bfnm{David~L}\binits{D.~L.}}
(\byear{1995}).
\btitle{Examples of basis pursuit}.
In \bbooktitle{SPIE's 1995 International Symposium on Optical Science,
  Engineering, and Instrumentation}
\bpages{564--574}.
\bpublisher{International Society for Optics and Photonics}.
\end{binproceedings}
\endbibitem

\bibitem{efron2004least}
\begin{barticle}[author]
\bauthor{\bsnm{Efron},~\bfnm{Bradley}\binits{B.}},
  \bauthor{\bsnm{Hastie},~\bfnm{Trevor}\binits{T.}},
  \bauthor{\bsnm{Johnstone},~\bfnm{Iain}\binits{I.}},
  \bauthor{\bsnm{Tibshirani},~\bfnm{Robert}\binits{R.}} \betal{et~al.}
(\byear{2004}).
\btitle{Least angle regression}.
\bjournal{The Annals of statistics}
\bvolume{32}
\bpages{407--499}.
\end{barticle}
\endbibitem

\bibitem{fan2001variable}
\begin{barticle}[author]
\bauthor{\bsnm{Fan},~\bfnm{Jianqing}\binits{J.}} \AND
  \bauthor{\bsnm{Li},~\bfnm{Runze}\binits{R.}}
(\byear{2001}).
\btitle{Variable selection via nonconcave penalized likelihood and its oracle
  properties}.
\bjournal{Journal of the American statistical Association}
\bvolume{96}
\bpages{1348--1360}.
\end{barticle}
\endbibitem

\bibitem{fan2004nonconcave}
\begin{barticle}[author]
\bauthor{\bsnm{Fan},~\bfnm{Jianqing}\binits{J.}} \AND
  \bauthor{\bsnm{Peng},~\bfnm{Heng}\binits{H.}}
(\byear{2004}).
\btitle{Nonconcave penalized likelihood with a diverging number of parameters}.
\bjournal{The Annals of Statistics}
\bvolume{32}
\bpages{928--961}.
\end{barticle}
\endbibitem

\bibitem{fan2014strong}
\begin{barticle}[author]
\bauthor{\bsnm{Fan},~\bfnm{Jianqing}\binits{J.}},
  \bauthor{\bsnm{Xue},~\bfnm{Lingzhou}\binits{L.}} \AND
  \bauthor{\bsnm{Zou},~\bfnm{Hui}\binits{H.}}
(\byear{2014}).
\btitle{Strong oracle optimality of folded concave penalized estimation}.
\bjournal{Annals of statistics}
\bvolume{42}
\bpages{819}.
\end{barticle}
\endbibitem

\bibitem{hartman1959functions}
\begin{barticle}[author]
\bauthor{\bsnm{Hartman},~\bfnm{Philip}\binits{P.}}
(\byear{1959}).
\btitle{On functions representable as a difference of convex functions}.
\bjournal{Pacific Journal of Mathematics}
\bvolume{9}
\bpages{707--713}.
\end{barticle}
\endbibitem

\bibitem{hiriart1985generalized}
\begin{bincollection}[author]
\bauthor{\bsnm{Hiriart-Urruty},~\bfnm{J-B}\binits{J.-B.}}
(\byear{1985}).
\btitle{Generalized Differentiability/Duality and Optimization for Problems
  Dealing with Differences of Convex Functions}.
In \bbooktitle{Convexity and duality in optimization}
\bpages{37--70}.
\bpublisher{Springer}.
\end{bincollection}
\endbibitem

\bibitem{horst1999dc}
\begin{barticle}[author]
\bauthor{\bsnm{Horst},~\bfnm{Reiner}\binits{R.}} \AND
  \bauthor{\bsnm{Thoai},~\bfnm{Nguyen~V}\binits{N.~V.}}
(\byear{1999}).
\btitle{DC programming: overview}.
\bjournal{Journal of Optimization Theory and Applications}
\bvolume{103}
\bpages{1--43}.
\end{barticle}
\endbibitem

\bibitem{huang2008adaptive}
\begin{barticle}[author]
\bauthor{\bsnm{Huang},~\bfnm{Jian}\binits{J.}},
  \bauthor{\bsnm{Ma},~\bfnm{Shuangge}\binits{S.}} \AND
  \bauthor{\bsnm{Zhang},~\bfnm{Cun-Hui}\binits{C.-H.}}
(\byear{2008}).
\btitle{Adaptive Lasso for sparse high-dimensional regression models}.
\bjournal{Statistica Sinica}
\bpages{1603--1618}.
\end{barticle}
\endbibitem

\bibitem{javanmard2014hypothesis}
\begin{barticle}[author]
\bauthor{\bsnm{Javanmard},~\bfnm{Adel}\binits{A.}} \AND
  \bauthor{\bsnm{Montanari},~\bfnm{Andrea}\binits{A.}}
(\byear{2014}).
\btitle{Hypothesis testing in high-dimensional regression under the gaussian
  random design model: Asymptotic theory}.
\bjournal{IEEE Transactions on Information Theory}
\bvolume{60}
\bpages{6522--6554}.
\end{barticle}
\endbibitem

\bibitem{javanmard2014confidence}
\begin{barticle}[author]
\bauthor{\bsnm{Javanmard},~\bfnm{Adel}\binits{A.}} \AND
  \bauthor{\bsnm{Montanari},~\bfnm{Andrea}\binits{A.}}
(\byear{2014}).
\btitle{Confidence intervals and hypothesis testing for high-dimensional
  regression.}
\bjournal{Journal of Machine Learning Research}
\bvolume{15}
\bpages{2869--2909}.
\end{barticle}
\endbibitem

\bibitem{lehmann1998theory}
\begin{barticle}[author]
\bauthor{\bsnm{Lehmann},~\bfnm{EL}\binits{E.}} \AND
  \bauthor{\bsnm{Casella},~\bfnm{G}\binits{G.}}
(\byear{1998}).
\btitle{Theory of Point Estimation, Springer-Verlag}.
\bjournal{New York}.
\end{barticle}
\endbibitem

\bibitem{loh2013regularized}
\begin{binproceedings}[author]
\bauthor{\bsnm{Loh},~\bfnm{Po-Ling}\binits{P.-L.}} \AND
  \bauthor{\bsnm{Wainwright},~\bfnm{Martin~J}\binits{M.~J.}}
(\byear{2013}).
\btitle{Regularized M-estimators with nonconvexity: Statistical and algorithmic
  theory for local optima}.
In \bbooktitle{Advances in Neural Information Processing Systems}
\bpages{476--484}.
\end{binproceedings}
\endbibitem

\bibitem{lv2009unified}
\begin{barticle}[author]
\bauthor{\bsnm{Lv},~\bfnm{Jinchi}\binits{J.}} \AND
  \bauthor{\bsnm{Fan},~\bfnm{Yingying}\binits{Y.}}
(\byear{2009}).
\btitle{A unified approach to model selection and sparse recovery using
  regularized least squares}.
\bjournal{The Annals of Statistics}
\bpages{3498--3528}.
\end{barticle}
\endbibitem

\bibitem{mazumder2011sparsenet}
\begin{barticle}[author]
\bauthor{\bsnm{Mazumder},~\bfnm{Rahul}\binits{R.}},
  \bauthor{\bsnm{Friedman},~\bfnm{Jerome~H}\binits{J.~H.}} \AND
  \bauthor{\bsnm{Hastie},~\bfnm{Trevor}\binits{T.}}
(\byear{2011}).
\btitle{Sparsenet: Coordinate descent with nonconvex penalties}.
\bjournal{Journal of the American Statistical Association}
\bvolume{106}
\bpages{1125--1138}.
\end{barticle}
\endbibitem

\bibitem{nouiehed2017pervasiveness}
\begin{barticle}[author]
\bauthor{\bsnm{Nouiehed},~\bfnm{Maher}\binits{M.}},
  \bauthor{\bsnm{Pang},~\bfnm{Jong-Shi}\binits{J.-S.}} \AND
  \bauthor{\bsnm{Razaviyayn},~\bfnm{Meisam}\binits{M.}}
(\byear{2017}).
\btitle{On the Pervasiveness of Difference-Convexity in Optimization and
  Statistics}.
\bjournal{arXiv preprint arXiv:1704.03535}.
\end{barticle}
\endbibitem

\bibitem{pang2016computing}
\begin{barticle}[author]
\bauthor{\bsnm{Pang},~\bfnm{Jong-Shi}\binits{J.-S.}},
  \bauthor{\bsnm{Razaviyayn},~\bfnm{Meisam}\binits{M.}} \AND
  \bauthor{\bsnm{Alvarado},~\bfnm{Alberth}\binits{A.}}
(\byear{2016}).
\btitle{Computing B-stationary points of nonsmooth DC programs}.
\bjournal{Mathematics of Operations Research}
\bvolume{42}
\bpages{95--118}.
\end{barticle}
\endbibitem

\bibitem{rockafellar2015convex}
\begin{bbook}[author]
\bauthor{\bsnm{Rockafellar},~\bfnm{Ralph~Tyrell}\binits{R.~T.}}
(\byear{2015}).
\btitle{Convex analysis}.
\bpublisher{Princeton university press}.
\end{bbook}
\endbibitem

\bibitem{sriperumbudur2012proof}
\begin{barticle}[author]
\bauthor{\bsnm{Sriperumbudur},~\bfnm{Bharath~K}\binits{B.~K.}} \AND
  \bauthor{\bsnm{Lanckriet},~\bfnm{Gert~RG}\binits{G.~R.}}
(\byear{2012}).
\btitle{A proof of convergence of the concave-convex procedure using zangwill's
  theory}.
\bjournal{Neural computation}
\bvolume{24}
\bpages{1391--1407}.
\end{barticle}
\endbibitem

\bibitem{tao1997convex}
\begin{barticle}[author]
\bauthor{\bsnm{Tao},~\bfnm{Pham~Dinh}\binits{P.~D.}} \AND
  \bauthor{\bsnm{An},~\bfnm{Le~Thi~Hoai}\binits{L.~T.~H.}}
(\byear{1997}).
\btitle{Convex analysis approach to dc programming: Theory, algorithms and
  applications}.
\bjournal{Acta Mathematica Vietnamica}
\bvolume{22}
\bpages{289--355}.
\end{barticle}
\endbibitem

\bibitem{tao2005dc}
\begin{barticle}[author]
\bauthor{\bsnm{Tao},~\bfnm{Pham~Dinh}\binits{P.~D.}} \betal{et~al.}
(\byear{2005}).
\btitle{The DC (difference of convex functions) programming and DCA revisited
  with DC models of real world nonconvex optimization problems}.
\bjournal{Annals of operations research}
\bvolume{133}
\bpages{23--46}.
\end{barticle}
\endbibitem

\bibitem{tibshirani1996regression}
\begin{barticle}[author]
\bauthor{\bsnm{Tibshirani},~\bfnm{Robert}\binits{R.}}
(\byear{1996}).
\btitle{Regression shrinkage and selection via the lasso}.
\bjournal{Journal of the Royal Statistical Society. Series B (Methodological)}
\bpages{267--288}.
\end{barticle}
\endbibitem

\bibitem{tuy1987global}
\begin{barticle}[author]
\bauthor{\bsnm{Tuy},~\bfnm{Hoang}\binits{H.}}
(\byear{1987}).
\btitle{Global minimization of a difference of two convex functions}.
\bjournal{Nonlinear Analysis and Optimization}
\bpages{150--182}.
\end{barticle}
\endbibitem

\bibitem{van2014asymptotically}
\begin{barticle}[author]
\bauthor{\bparticle{Van~de} \bsnm{Geer},~\bfnm{Sara}\binits{S.}},
  \bauthor{\bsnm{B{\"u}hlmann},~\bfnm{Peter}\binits{P.}},
  \bauthor{\bsnm{Ritov},~\bfnm{Ya'acov}\binits{Y.}} \AND
  \bauthor{\bsnm{Dezeure},~\bfnm{Ruben}\binits{R.}}
(\byear{2014}).
\btitle{On asymptotically optimal confidence regions and tests for
  high-dimensional models}.
\bjournal{The Annals of Statistics}
\bvolume{42}
\bpages{1166--1202}.
\end{barticle}
\endbibitem

\bibitem{wainwright2009sharp}
\begin{barticle}[author]
\bauthor{\bsnm{Wainwright},~\bfnm{Martin~J}\binits{M.~J.}}
(\byear{2009}).
\btitle{Sharp thresholds for High-Dimensional and noisy sparsity recovery using
  $l_1$-Constrained Quadratic Programming (Lasso)}.
\bjournal{IEEE transactions on information theory}
\bvolume{55}
\bpages{2183--2202}.
\end{barticle}
\endbibitem

\bibitem{wang2013calibrating}
\begin{barticle}[author]
\bauthor{\bsnm{Wang},~\bfnm{Lan}\binits{L.}},
  \bauthor{\bsnm{Kim},~\bfnm{Yongdai}\binits{Y.}} \AND
  \bauthor{\bsnm{Li},~\bfnm{Runze}\binits{R.}}
(\byear{2013}).
\btitle{Calibrating non-convex penalized regression in ultra-high dimension}.
\bjournal{Annals of statistics}
\bvolume{41}
\bpages{2505}.
\end{barticle}
\endbibitem

\bibitem{wang2014optimal}
\begin{barticle}[author]
\bauthor{\bsnm{Wang},~\bfnm{Zhaoran}\binits{Z.}},
  \bauthor{\bsnm{Liu},~\bfnm{Han}\binits{H.}} \AND
  \bauthor{\bsnm{Zhang},~\bfnm{Tong}\binits{T.}}
(\byear{2014}).
\btitle{Optimal computational and statistical rates of convergence for sparse
  nonconvex learning problems}.
\bjournal{Annals of statistics}
\bvolume{42}
\bpages{2164}.
\end{barticle}
\endbibitem

\bibitem{yuille2003concave}
\begin{barticle}[author]
\bauthor{\bsnm{Yuille},~\bfnm{Alan~L}\binits{A.~L.}} \AND
  \bauthor{\bsnm{Rangarajan},~\bfnm{Anand}\binits{A.}}
(\byear{2003}).
\btitle{The concave-convex procedure}.
\bjournal{Neural computation}
\bvolume{15}
\bpages{915--936}.
\end{barticle}
\endbibitem

\bibitem{zhang2010nearly}
\begin{barticle}[author]
\bauthor{\bsnm{Zhang},~\bfnm{Cun-Hui}\binits{C.-H.}}
(\byear{2010}).
\btitle{Nearly unbiased variable selection under minimax concave penalty}.
\bjournal{The Annals of statistics}
\bvolume{38}
\bpages{894--942}.
\end{barticle}
\endbibitem

\bibitem{zhang2014confidence}
\begin{barticle}[author]
\bauthor{\bsnm{Zhang},~\bfnm{Cun-Hui}\binits{C.-H.}} \AND
  \bauthor{\bsnm{Zhang},~\bfnm{Stephanie~S}\binits{S.~S.}}
(\byear{2014}).
\btitle{Confidence intervals for low dimensional parameters in high dimensional
  linear models}.
\bjournal{Journal of the Royal Statistical Society: Series B (Statistical
  Methodology)}
\bvolume{76}
\bpages{217--242}.
\end{barticle}
\endbibitem

\bibitem{zhang2014minimization}
\begin{barticle}[author]
\bauthor{\bsnm{Zhang},~\bfnm{Shuai}\binits{S.}} \AND
  \bauthor{\bsnm{Xin},~\bfnm{Jack}\binits{J.}}
(\byear{2014}).
\btitle{Minimization of Transformed L\_1 Penalty: Theory, Difference of Convex
  Function Algorithm, and Robust Application in Compressed Sensing}.
\bjournal{arXiv preprint arXiv:1411.5735}.
\end{barticle}
\endbibitem

\bibitem{zhang2010analysis}
\begin{barticle}[author]
\bauthor{\bsnm{Zhang},~\bfnm{Tong}\binits{T.}}
(\byear{2010}).
\btitle{Analysis of multi-stage convex relaxation for sparse regularization}.
\bjournal{Journal of Machine Learning Research}
\bvolume{11}
\bpages{1081--1107}.
\end{barticle}
\endbibitem

\bibitem{zhang2013multi}
\begin{barticle}[author]
\bauthor{\bsnm{Zhang},~\bfnm{Tong}\binits{T.}}
(\byear{2013}).
\btitle{Multi-stage convex relaxation for feature selection}.
\bjournal{Bernoulli}
\bvolume{19}
\bpages{2277--2293}.
\end{barticle}
\endbibitem

\bibitem{zhao2006model}
\begin{barticle}[author]
\bauthor{\bsnm{Zhao},~\bfnm{Peng}\binits{P.}} \AND
  \bauthor{\bsnm{Yu},~\bfnm{Bin}\binits{B.}}
(\byear{2006}).
\btitle{On model selection consistency of Lasso}.
\bjournal{Journal of Machine learning research}
\bvolume{7}
\bpages{2541--2563}.
\end{barticle}
\endbibitem

\bibitem{zou2006adaptive}
\begin{barticle}[author]
\bauthor{\bsnm{Zou},~\bfnm{Hui}\binits{H.}}
(\byear{2006}).
\btitle{The adaptive lasso and its oracle properties}.
\bjournal{Journal of the American statistical association}
\bvolume{101}
\bpages{1418--1429}.
\end{barticle}
\endbibitem

\bibitem{zou2008one}
\begin{barticle}[author]
\bauthor{\bsnm{Zou},~\bfnm{Hui}\binits{H.}} \AND
  \bauthor{\bsnm{Li},~\bfnm{Runze}\binits{R.}}
(\byear{2008}).
\btitle{One-step sparse estimates in nonconcave penalized likelihood models}.
\bjournal{Annals of statistics}
\bvolume{36}
\bpages{1509}.
\end{barticle}
\endbibitem

\end{thebibliography}

\clearpage
\appendix
\section{Properties of DC programming}\label{app:dcproperty}

The following are some known properties of the DC functions \citep{tao1997convex, horst1999dc}.
\begin{enumerate}
\item Every DC function has a nonnegative DC decomposition; that is for a DC function $f$, there exists a decomposition, $f = g - h$, where both $g$ and $h$ are nonnegative and convex.
\item Every $C^1$ (i.e., functions with continuously first order derivatives) function with a Lipschitz gradient is a DC function.
\item Every twice continuously differentiable function is a DC function.
\item Every continuous function on a convex set is a limit of a sequence of uniformly converging DC functions.
\item Let $f_i$ be DC functions for $i = 1, \cdots, m$.
The DC functions are {\it closed} under the following operations:
\begin{itemize}
\item summation: $\sum_{i = 1}^m\lambda_if_i(x)$, for $\lambda_i \in \mathbb{R}$, $i = 1, \cdots, m$
\item maximization: $\max_{i = 1, \cdots, m}f_i(x)$
\item minimization: $\min_{i = 1, \cdots, m}f_i(x)$
\item product: $\prod_{i = 1, \cdots, m}f_i(x)$
\end{itemize}

\item A locally DC function that is defined in $\mathbb{R}^n$ is a DC function.
\item The following statements about a DC program are equivalent:
\begin{itemize}
\item $\sup\{f(x): x\in C\}$, function $f$ and set $C$ are convex
\item $\inf\{g(x) - h(x): x \in \mathbb{R}^n\}$, functions $g$ and $h$ are convex
\item $\inf\{g(x) - h(x): x\in C, f_1(x) - f_1(x) \leq 0\}$, functions $g$, $h$, $f_1$, and $f_2$ and set $C$ are all convex
\end{itemize}

\end{enumerate}

Regarding the optimal solutions in the DC programming, the following have been developed in the literature \citep{tao1997convex, horst1999dc}.
\begin{definition}
{\bf ($\epsilon$-subdifferential)}
For a convex function $g(x)$, and $\epsilon > 0$, the $\epsilon$-subdifferential of function $g(x)$ at point $x_0$ is denoted by $\partial_\epsilon g(x_0)$ and is defined as follows:
$$
\partial_\epsilon g(x_0)= \{\nu \in \mathbb{R}^n | g(x) \geq g(x_0) + \langle x- x_0, \nu \rangle- \epsilon\}.
$$
\end{definition}
One can verify that the subgradient \cite[Chapter 23]{rockafellar2015convex} (which is denoted by $\partial g(x_0)$) of function $g(x)$ at $x_0$ is the $0$-subdifferential (i.e., $\epsilon = 0$).
\begin{theorem}
\label{thm:globalopt}
{\bf (Global optimality condition)} A point $x^\ast$ is a global optimal if and only if (iff) $\partial_\epsilon h(x^\ast) \subset \partial_\epsilon g(x^\ast)$ for any $\epsilon > 0$.
\end{theorem}
\begin{theorem}
\label{thm:localopt}
{\bf (Local optimality condition)} A point $x^\ast$ is a local optimal if $\partial h(x^\ast) \subset \text{int }\partial g(x^\ast)$, where $\text{int }\partial g(x^\ast)$ represents the interior of the set $\partial g(x^\ast)$.
\end{theorem}

\section{Proofs in Section \ref{MAINRESULTS} }\label{appA}
\subsection{Proof of Theorem \ref{thm:estimation}}\label{proof:estimation}
\begin{proof}
Since $\hat{\beta}^{\lambda}_{n}$ is d-stationary, the directional derivatives should be nonnegative in all directions, especially in the direction of $\beta^{\ast} -\hat{\beta}^{\lambda}_n$:
\[-\frac{1}{n}(\beta^{\ast} -\hat{\beta}^{\lambda}_n)^TX^T(Y - X\hat{\beta}^{\lambda}_n) + P'_{\lambda}(\hat{\beta}^{\lambda}_n; \beta^{\ast} -\hat{\beta}^{\lambda}_n) \geq 0,\]
where $P'_{\lambda}(\hat{\beta}^{\lambda}_n; \beta^{\ast} -\hat{\beta}^{\lambda}_n)$ is the directional derivative for the penalty function $P_{\lambda} = \lambda \|\beta\|_1 - h_{\lambda}(\beta)$ at $\hat{\beta}^{\lambda}_n$ in the direction of $\beta^{\ast} -\hat{\beta}^{\lambda}_n$.

According to Lemma \ref{FOClemma}, there exists a subgradient $z \in \partial g(\hat{\beta}^{\lambda}_n)$, where  $\partial g(\hat{\beta}^{\lambda}_n)$ is the set of subgradient of $g(\beta) = \|\beta\|_1$ at $\hat{\beta}^{\lambda}_n$, such that:
\begin{equation}
 \nabla L(\hat{\beta}^{\lambda}_n) +  \lambda z -  \nabla h_\lambda(\hat{\beta}^{\lambda}_n) = 0,
\end{equation}

\noindent Multiplying by $(\beta^{\ast} -\hat{\beta}^{\lambda}_n)^T$ on both side and plugging in $Y = X\beta^{\ast} + \epsilon$, we have
\[-\frac{1}{n}(\beta^{\ast} -\hat{\beta}^{\lambda}_n)^TX^TX(\beta^{\ast} -\hat{\beta}^{\lambda}_n) - \frac{1}{n}(\beta^{\ast} -\hat{\beta}^{\lambda}_n)^TX^T\epsilon + P'_{\lambda}(\hat{\beta}^{\lambda}_n; \beta^{\ast} -\hat{\beta}^{\lambda}_n) = 0,\]
where without ambiguity, we let $P'_{\lambda}(\hat{\beta}^{\lambda}_n; \beta^{\ast} -\hat{\beta}^{\lambda}_n) = (\beta^{\ast} -\hat{\beta}^{\lambda}_n)^T(\lambda z -  \nabla h_\lambda(\hat{\beta}^{\lambda}_n))$ since the true directional derivative for the penalty $P_{\lambda} = \lambda \|\beta\|_1 - h_{\lambda}(\beta)$ at $\hat{\beta}^{\lambda}_n$ in the direction of $\beta^{\ast} -\hat{\beta}^{\lambda}_n$ is greater than $(\beta^{\ast} -\hat{\beta}^{\lambda}_n)^T(\lambda z -  \nabla h_\lambda(\hat{\beta}^{\lambda}_n))$. We thus will have
\begin{equation}
\label{e001}
\begin{split}
\frac{1}{n}(\beta^{\ast} -\hat{\beta}^{\lambda}_n)^TX^TX(\beta^{\ast}  - \hat{\beta}^{\lambda}_n) = - \frac{1}{n}(\beta^{\ast} -\hat{\beta}^{\lambda}_n)^TX^T\epsilon + P'_{\lambda}(\hat{\beta}^{\lambda}_n; \beta^{\ast} -\hat{\beta}^{\lambda}_n),
\end{split}
\end{equation}
which implies we have the following hold for $j \notin S$

\begin{equation}
\begin{split}
&-\frac{1}{n}|(\hat{\beta}^{\lambda}_n)_j|\text{sign}((\hat{\beta}^{\lambda}_n)_j)X_j^TX(\beta^{\ast}  - \hat{\beta}^{\lambda}_n)\\
=& \frac{1}{n}|(\hat{\beta}^{\lambda}_n)_j|\text{sign}((\hat{\beta}^{\lambda}_n)_j)X_j^T\epsilon - \lambda |(\hat{\beta}^{\lambda}_n)_j| + h'_\lambda((\hat{\beta}^{\lambda}_n)_j)|(\hat{\beta}^{\lambda}_n)_j|\text{sign}((\hat{\beta}^{\lambda}_n)_j)\\
\leq& -c\lambda |(\hat{\beta}^{\lambda}_n)_j|,
\end{split}
\end{equation}
where the ``$\leq$'' follows from Assumption \ref{nonnegativity}.

For $j \in S$

\begin{equation}
\begin{split}
&\frac{1}{n}(\beta^{\ast} -\hat{\beta}^{\lambda}_n)_jX_j^TX(\beta^{\ast}  - \hat{\beta}^{\lambda}_n)\\
=& - \frac{1}{n}(\beta^{\ast} -\hat{\beta}^{\lambda}_n)_jX_j^T\epsilon
+ \lambda (\beta^{\ast} -\hat{\beta}^{\lambda}_n)_j z_j
- h'_\lambda((\hat{\beta}^{\lambda}_n)_j)(\beta^{\ast} -\hat{\beta}^{\lambda}_n)_j\\
\leq& \frac{5}{2}\lambda |(\beta^{\ast} -\hat{\beta}^{\lambda}_n)_j|,
\end{split}
\end{equation}
where the ``$\leq$'' follows from $\lambda \geq \frac{2\|X^T\epsilon\|_{\infty}}{n} $ and Assumption \ref{hderivative}.

Let $\nu = \beta^{\ast} - \hat{\beta}^{\lambda}_{n}$, we thus will have
\begin{equation}
\begin{split}
\frac{1}{n}(\beta^{\ast} -\hat{\beta}^{\lambda}_n)^TX^TX(\beta^{\ast}  - \hat{\beta}^{\lambda}_n)
\leq  - c\lambda\|(\nu)_{S^c}\|_1 + \frac{5}{2}\lambda\|(\nu)_{S}\|_1,
\end{split}
\end{equation}
where the inequality follows from Assumption \ref{hderivative}.

Since the left hand side of the above is nonnegative, we will have $\nu = \beta^{\ast} - \hat{\beta}^{\lambda}_{n} \in \mathcal{C}$.
Under the restricted strong convexity condition, we will have
\begin{equation}
\label{e01}
\begin{split}
\gamma \|\nu\|_2^2 \leq \frac{1}{n}(\beta^{\ast} -\hat{\beta}^{\lambda}_n)^TX^TX(\beta^{\ast}  - \hat{\beta}^{\lambda}_n) &\leq \frac{5}{2}\lambda\|(\nu)_{S}\|_1
\end{split}
\end{equation}
Thus we will further have
\[\|\nu\|_2^2 \leq \frac{5}{2\gamma}\lambda\|(\nu)_{S}\|_1 \leq  \frac{5}{2\gamma}\lambda\sqrt{|S|}\|(\nu)\|_2 ,\]
from which we will have the upper bound
\begin{equation}
\begin{split}
 \|\nu\|_2 \leq  \frac{5}{2\gamma}\lambda\sqrt{|S|}  \propto \frac{\lambda \sqrt{|S|}}{\gamma}
 \end{split}
\end{equation}

\end{proof}
\subsection{Proof of Corollary \ref{cor:estimation}}\label{proof:corestimation}
\begin{proof}
From the proof of Theorem \ref{thm:estimation}, let $\nu = \beta^{\ast} - \hat{\beta}^{\lambda}_{n}$,we have
\begin{equation}
\begin{split}
&\|\frac{1}{\sqrt{n}}X(\beta^{\ast}  - \hat{\beta}^{\lambda}_n)\|^2_2 \\
 = &\frac{1}{n}(\beta^{\ast} -\hat{\beta}^{\lambda}_n)^TX^TX(\beta^{\ast}  - \hat{\beta}^{\lambda}_n) \\
\leq&  - c\lambda\|(\nu)_{S^c}\|_1 + \frac{5}{2}\lambda\|(\nu)_{S}\|_1\\
\leq &(\frac{5}{2}\lambda)^2 \frac{1}{\gamma}\sqrt{|S|}
\end{split}
\end{equation}
\end{proof}
\subsection{Proof of Corollary \ref{cor:Bon}}\label{proof:corBon}
\begin{proof}
Since $\epsilon_i$ for $i = 1, \cdots, n$ are from sub-Gaussian distribution with parameter $\sigma^2$,
\[
\mathbb{P}\left(\frac{x_j^T\epsilon}{n} \geq t \right)
\leq 2\exp{(-\frac{nt^2}{2\sigma^2})}\]
By Bonferroni bound, we will have
\[\mathbb{P}(\frac{\|X^T\epsilon\|_{\infty}}{n} \geq t) \leq 2\exp{(-\frac{nt^2}{2\sigma^2} + \log{p})}\]
By setting $t = \sigma\sqrt{\frac{\tau \log{p}}{n}}$ for some $\tau \geq 2$, we will be able to have
\[\mathbb{P}(\frac{\|X^T\epsilon\|_{\infty}}{n} \geq t) \leq  2\exp{(-\frac{\tau - 2}{2}\log{p})}.\]
Thus:
\begin{equation}
\begin{split}
\|\hat{\beta}^{\lambda}_n - \beta^{\ast}\|_2 \lesssim \frac{5}{\gamma}\sigma\sqrt{\frac{\tau |S|\log{p}}{n}}
\end{split}
\end{equation}
\end{proof}

\subsection{Proof of Lemma \ref{FOClemma}}\label{proof:FOClemma}
\begin{proof}
We will first prove the necessity. Since $\beta_0$ is a d-stationary solution to the objective function $F(\beta)$, we will have
\[F^{\prime}(\beta_0,\beta - \beta_0) \geq 0,\]
for any $\beta \in \mathbb{R}^p$, where $F^{\prime}(\beta_0,\beta - \beta_0)$ denotes the directional derivative in the direction of $\beta - \beta_0$. For any $i = 1, \cdots, p$, let $\beta^{i+} = \beta_0 + e_i$, where $e_i \in \mathbb{R}^p$ denotes the unit vector with $1$ in the ith position and $0$ everywhere else. Let $\beta^{i-} = \beta_0 - e_i$. We will have:
\[F^{\prime}(\beta_0,\beta^{i+} - \beta_0) \geq 0,\]
\[F^{\prime}(\beta_0,\beta^{i-} - \beta_0) \geq 0,\]
which implies
\begin{itemize}
\item For $i$ such that $\beta_{0i} \neq 0$,
\[\nabla L(\beta_0)_i - \nabla h(\beta_0)_i + \text{sign}(\beta_{0i}) \geq 0,\]
\[-\nabla L(\beta_0)_i + \nabla h(\beta_0)_i -  \text{sign}(\beta_{0i}) \geq 0,\]
where $\text{sign}(x) = 1$ if $x > 0$, $\text{sign}(x) = -1$ if $x < 0$.
\item For $i$ such that $\beta_{0i} = 0$,
\[\nabla L(\beta_0)_i - \nabla h(\beta_0)_i + 1\geq 0,\]
\[-\nabla L(\beta_0)_i + \nabla h(\beta_0)_i + 1\geq 0.\]
\end{itemize}
 We thus conclude that

 \begin{itemize}
\item For $i$ such that $\beta_{0i} \neq 0$,
\[\nabla L(\beta_0)_i - \nabla h(\beta_0)_i + \text{sign}(\beta_{0i}) = 0,\]
\item For $i$ such that $\beta_{0i} = 0$,
\[|\nabla L(\beta_0)_i - \nabla h(\beta_0)_i| \leq 1.\]
\end{itemize}
Thus for $z$ such that $z_i = \text{sign}(\beta_{0i})$ when $\beta_{0i} \neq 0$ and $z_i = -(\nabla L(\beta_0)_i - \nabla h(\beta_0)_i)$ when $\beta_{0i} = 0$ is what we need.

On the other hand, if there exists some $z \in \partial g(\beta_0)$, where  $\partial g(\beta_0)$ is the set of subgradient of $g(\beta)$ at $\beta_0$, such that:
\begin{equation}
 \nabla L(\beta_0) +  z -  \nabla h(\beta_0) = 0,
\end{equation}
 \begin{itemize}
\item For $i$ such that $\beta_{0i} \neq 0$,
\[ z_i = \text{sign}(\beta_{0i}) = 1\text{ or } -1,\]
\item For $i$ such that $\beta_{0i} = 0$,
\[-1 \leq z = -(\nabla L(\beta_0)_i - \nabla h(\beta_0)_i) \leq 1.\]
\end{itemize}

 \begin{itemize}
\item For $i$ such that $\beta_{0i} \neq 0$,
\[ F^{\prime}(\beta_0,\beta^{i+} - \beta_0) = 0,\]
\[ F^{\prime}(\beta_0,\beta^{i-} - \beta_0) = 0.\]
\item For $i$ such that $\beta_{0i} = 0$,
\begin{equation}
\begin{split}
F^{\prime}(\beta_0,\beta^{i+} - \beta_0) &= \nabla L(\beta_0)_i - \nabla h(\beta_0)_i + 1\\
&= \nabla L(\beta_0)_i - \nabla h(\beta_0)_i +z + 1 - z\\
&= 0 + 1 - z \geq 0
\end{split}
\end{equation}

\begin{equation}
\begin{split}
F^{\prime}(\beta_0,\beta^{i-} - \beta_0) &= -\nabla L(\beta_0)_i + \nabla h(\beta_0)_i + 1\\
&= -\nabla L(\beta_0)_i + \nabla h(\beta_0)_i - z + 1 + z\\
&= 0 + 1 + z \geq 0
\end{split}
\end{equation}
\end{itemize}
We thus conclude that the directional derivative of the objective function at $\beta_0$ is always nonnegative in any direction.
This complete the proof.
\end{proof}

\begin{remark}
From the proof Lemma \ref{FOClemma}, we can derive similar conditions for ``local maximals'' for $\tilde{\beta}$ satisfying the following:
\begin{equation}
F^{\prime}(\tilde{\beta},\beta - \tilde{\beta}) \leq 0.
\end{equation}

\begin{itemize}
\item For $i$ such that $\tilde{\beta}_{i} \neq 0$,
\[ F^{\prime}(\tilde{\beta},\beta^{i+} - \tilde{\beta}) = 0,\]
\[ F^{\prime}(\tilde{\beta},\beta^{i-} - \tilde{\beta}) = 0.\]
\item For $i$ such that $\tilde{\beta}_{i} = 0$,
\[\nabla L(\tilde{\beta})_i - \nabla h(\tilde{\beta})_i + 1\leq 0,\]
\[-\nabla L(\tilde{\beta})_i + \nabla h(\tilde{\beta})_i + 1\leq 0.\]
\end{itemize}
which implies that if the stationary solution $\tilde{\beta}$ to the FOC satisfies: $\min_{i = 1}^p\{\tilde{\beta}_{i}\} = 0$, it will only satisfy the condition for ``local'' minimals and thus be a d-stationary solution.
\end{remark}

\subsection{Proof of Lemma \ref{lem:cvxity}}\label{proof:lemcvxity}
\begin{proof}
Since $L(\beta) = \frac{1}{2n}\|Y - X\beta\|^2_2$ is quadratic and convex, we have
\[L(\beta_2) = L(\beta_1) + \nabla L(\beta_1)^T(\beta_2 - \beta_1) + \frac{1}{2}(\beta_2 - \beta_1)^T\nabla^2L(\beta_1)(\beta_2 - \beta_1),
\]
where $\nabla^2L(\beta_1)$ is the Hessian matrix of $L(\beta)$ at $\beta_1$. Since $\nu = \beta_1 - \beta_2 \in \mathcal{C}$ and Assumption \ref{RSC} holds on $\mathcal{C}$, we will further have
\[L(\beta_2) \geq L(\beta_1) + \nabla L(\beta_1)^T(\beta_2 - \beta_1) + \frac{\gamma}{2}\|\beta_2 - \beta_1\|^2_2.
\]
On the other hand, $h_{\lambda}(\beta)$ is convex with $0 \leq \eta^+ \leq \frac{h'_\lambda(t_2) - h'_\lambda(t_1)}{t_2 - t_1} \leq \eta^-$, we will have
\[h_{\lambda}(\beta_2) \leq h_{\lambda}(\beta_1) + \nabla h_{\lambda}(\beta_1)^T(\beta_2 - \beta_1) + \frac{\eta^-}{2}\|\beta_2 - \beta_1\|^2_2\]
By combining the above two inequalities, we will be able to get (\ref{cvxity}).
\end{proof}

\subsection{Proof of Theorem \ref{oracle}}\label{proof:oracle}
\begin{proof}
The first part is easy to see since the feasible region is convex and is a subset of $\mathcal{C}$, in which the strong convexity condition holds (Assumption \ref{RSC}) for the loss function (in our case, the least square loss function).
The minimizer to a strong problem is unique.

\noindent For the second conclusion, we first need to show that $X^T_SX_S$ is invertible and $\beta^O_S = (X^T_SX_S)^{-1}X^T_SY$. This follows easily from the Assumption \ref{RSC}, which implies $\gamma_{\min}(X^T_SX_S)$, the minimum eigenvalue of $X^T_SX_S$ is larger than $n\gamma$. We thus have $\beta^O_S = (X^T_SX_S)^{-1}X^T_SY$.
\begin{equation}
\begin{split}
\beta^O_S - \beta^\ast_S &= (X^T_SX_S)^{-1}X^T_SY - \beta^\ast\\
&= (X^T_SX_S)^{-1}X^T_S\epsilon
\end{split}
\end{equation}
Since $e_j(X^T_SX_S)^{-1}X^T_S\epsilon$, where $e_j \in \mathbb{R}^s$with all-zero elements except the $j$-th coordinate. Recall that $\epsilon$ has independent sub-Gaussian coordinates with the same variance parameter $\sigma^2$, we thus will have
\begin{equation}
\mathbb{P}\left(|e_j(X^T_SX_S)^{-1}X^T_S\epsilon| > t\right) \leq 2\exp{-\frac{t^2}{\|e_j(X^T_SX_S)^{-1}X^T_S\|^2_2\sigma^2}}.
\end{equation}
By using Bonferroni bound, the above implies
\begin{equation}
\mathbb{P}\left(\max_{j = 1,\cdots,s}|e_j(X^T_SX_S)^{-1}X^T_S\epsilon| > t\right) \leq 2s\exp{-\frac{t^2}{\|e_j(X^T_SX_S)^{-1}X^T_S\|^2_2\sigma^2}}.
\end{equation}
Taking $t = C\|e_j(X^T_SX_S)^{-1}X^T_S\epsilon\|_2\sigma\cdot \sqrt{2\log{s}}$ with $C > 0$, we will have
\begin{equation}
\begin{split}
\|\beta^O_S - \beta^\ast_S\|_\infty &= \|(X^T_SX_S)^{-1}X^T_S\epsilon\|_\infty \\
& = \max_{j = 1,\cdots,s}|e_j(X^T_SX_S)^{-1}X^T_S\epsilon| \\
&\leq C\|e_j(X^T_SX_S)^{-1}X^T_S\epsilon\|_2\sigma\cdot \sqrt{2\log{s}}
\end{split}
\end{equation}
hold with probability at least $1 - 2\exp{-C^2}/s$. Since for any $j \in \{1, \cdots, s\}$,

\begin{equation}
\begin{split}
\|e_j(X^T_SX_S)^{-1}X^T_S\epsilon\|_2^2
& = e_j(X^T_SX_S)^{-1}X^T_S\epsilon(e_j(X^T_SX_S)^{-1}X^T_S\epsilon^T)\\
& = e_j(X^T_SX_S)^{-1}e_j^T\\
&\leq 1/\gamma_{\min}(X^T_SX_S)\\
&\leq 1/n\gamma.
\end{split}
\end{equation}
This complete the proof since $\|\beta^O_{S^c} - \beta^\ast_{S^c}\|_\infty = 0$.
\end{proof}

\subsection{Proof of Lemma \ref{oracledstat}}\label{proof:oracledstat}
\begin{proof}
According to Theorem \ref{oracle}, we will have for $j \in S$,
\[|\beta^O_j| \geq |\beta^{\ast}_j| - \|\beta^O - \beta^{\ast}\|_{\infty} \geq \zeta,\]
which further implies $P'_j(\lambda,\beta^O_j) = 0$.

For $j \notin S$, we will have $h'_{\lambda}(\beta^O_j) = 0$ and since the errors are sub-Gaussian, there will exist $\xi^O_{S^c} \in  \partial\|\beta^O_{S^c}\|_1$ satisfying inequality (\ref{focOracle}) with high probability and $\|\xi^O_{S^c}\|_{\infty} \leq \frac{1}{10}c$ where $c$ is defined in the Assumption \ref{nonnegativity}.
\[\left(\nabla L_{n}(\beta^O)\right)_{S^c} + \lambda\xi^O_{S^c} = 0.\]
In order to see this, first we notice that $\beta^O_S = (X^T_SX_S)^{-1}X^T_SY$, $\beta^O_{S^c} = 0$. We thus will have $\left(\nabla L_{n}(\beta^O)\right)_{S^c} = - \frac{1}{n}X_{S^c}(Y - (X^T_SX_S)^{-1}X^T_SY)$. Plugging in the true model $ Y = X\beta^\ast + \epsilon$, we will have $\left(\nabla L_{n}(\beta^O)\right)_{S^c} = - \frac{1}{n}X_{S^c}(I - X_{S}(X^T_SX_S)^{-1}X^T_S)\epsilon$, where $(I - X_{S}(X^T_SX_S)^{-1}X^T_S)\epsilon$ is a vector of independent sub-Gaussian random variables. By using the Bonferroni bound, we will have the conclusion.
\end{proof}

\subsection{Proof of Lemma \ref{subset}}\label{proof:subset}
\begin{proof}
Since both $\hat{\beta}$ and $\beta^O$ are d-stationary, per Lemma \ref{oracledstat}, we have:
\[(\hat{\beta} - \beta^O)^T \left(\nabla f_{\lambda}(\beta^O) + \lambda \xi^O\right) \geq 0,\]
and $\hat{\beta}$ is d-stationary, there exists a $\hat{\xi} \in \partial\|\hat{\beta}\|_1$ ((where $\partial \|\hat{\beta}\|_1$ stands for the subgradient of function $\|\beta\|_1$ at $\beta = \hat{\beta}$)) such that
\[(\beta ^O- \hat{\beta})^T \left(\nabla f_{\lambda}(\hat{\beta}) + \lambda \hat{\xi}\right) \geq 0.\]
On the one hand,
\begin{equation}
\begin{split}
0 &\leq (\beta ^O- \hat{\beta})^T \left(\nabla f_{\lambda}(\hat{\beta}) + \lambda \hat{\xi}\right)\\
&\leq (\beta ^O- \hat{\beta})^T \left(\nabla f_{\lambda}(\hat{\beta})\right)
-  \lambda\|(\nu)_{S^c}\|_1
+ \lambda\|(\nu)_{S}\|_1\\
& = -\frac{1}{n}(\beta ^O- \hat{\beta})^TX^TX(\beta ^{\ast}  - \beta ^O
+ \beta ^O - \hat{\beta}) - \frac{1}{n}(\beta ^O- \hat{\beta})^TX^T\epsilon\\
 &- (\beta ^O- \hat{\beta})^T\nabla h_{\lambda}(\hat{\beta})-  \lambda\|(\nu)_{S^c}\|_1 + \lambda\|(\nu)_{S}\|_1\\
&\leq -\frac{1}{n}(\beta ^O- \hat{\beta})^TX^TX(\beta ^{\ast}  - \beta ^O
+ \beta ^O - \hat{\beta}) - \frac{1}{n}(\beta ^O- \hat{\beta})^TX^T\epsilon \\
 &+ \sum_{i \notin S}|\hat{\beta}_{i}||h'_{\lambda}(\hat{\beta}_{i})|-  \lambda\|(\nu)_{S^c}\|_1 + 2\lambda\|(\nu)_{S}\|_1.
\end{split}
\end{equation}
By rearranging the terms, we will have
\begin{equation}
\label{l1}
\begin{split}
&\frac{1}{n}(\beta ^O- \hat{\beta})^TX^TX(\beta ^{\ast}  - \beta ^O) + \frac{1}{n}(\beta ^O- \hat{\beta})^TX^TX(\beta ^O - \hat{\beta})
- \sum_{i \notin S}|\hat{\beta}_{i}||h'_{\lambda}(\hat{\beta}_{i}) |\\
  \leq & - \frac{1}{n}(\beta ^O- \hat{\beta})^TX^T\epsilon
  -  \lambda\|(\nu)_{S^c}\|_1 + 2\lambda\|(\nu)_{S}\|_1.
\end{split}
\end{equation}
On the other hand, according to the proof of Lemma \ref{oracledstat}, we will have
\[\left(\nabla h_{\lambda}(\beta^O))\right)_{S} = \lambda\text{sign}(\beta^O_S),\]
\[\left(\nabla h_{\lambda}(\beta^O))\right)_{S^c} = 0,\]
\[\lambda\xi^O_{S^c} = -\left(\nabla L_{n}(\beta^O)\right)_{S^c},\]
\[\|\xi^O_{S^c}\|_{\infty} \leq \frac{1}{2}.\]
By using the above facts, we will further obtain
\begin{equation}
\begin{split}
0 &\leq (\hat{\beta} - \beta^O)^T \left(\nabla f_{\lambda}(\beta^O) + \lambda \xi^O\right)\\
& = -\frac{1}{n}(\hat{\beta} - \beta^O)^TX^TX(\beta ^{\ast}  - \beta ^O) - \frac{1}{n}(\hat{\beta} - \beta^O)^TX^T\epsilon + (\hat{\beta} - \beta^O)^T_{S^c}\lambda\xi^O_{S^c}\\
&\leq -\frac{1}{n}(\hat{\beta} - \beta^O)^TX^TX(\beta ^{\ast}  - \beta ^O) - \frac{1}{n}(\hat{\beta} - \beta^O)^TX^T\epsilon + \frac{1}{10}c\lambda\|(\nu)_{S^c}\|_1.
\end{split}
\end{equation}
By rearranging the terms, we will have
\begin{equation}
\label{l2}
\begin{split}
\frac{1}{n}(\hat{\beta} - \beta^O)^TX^T\epsilon - \frac{1}{10}c\lambda\|(\nu)_{S^c}\|_1 \leq \frac{1}{n}(\beta ^O- \hat{\beta})^TX^TX(\beta ^{\ast}  - \beta ^O)
\end{split}
\end{equation}
Plugging inequality (\ref{l2}) to inequality (\ref{l1}), we will have
\begin{equation}
 \frac{1}{n}(\beta ^O- \hat{\beta})^TX^TX(\beta ^O - \hat{\beta})
 \leq  \sum_{i \notin S}|\hat{\beta}_{i}||h'_{\lambda}(\hat{\beta}_{i}) |
   -   \lambda\|(\nu)_{S^c}\|_1 + \frac{1}{10}c\lambda\|(\nu)_{S^c}\|_1+ 2\lambda\|(\nu)_{S}\|_1.
\end{equation}
Under the Assumption \ref{nonnegativity} with use of Bonferroni bound as in Corollary \ref{cor:Bon}, we will have
\begin{equation}
 0 \leq \frac{1}{n}(\beta ^O- \hat{\beta})^TX^TX(\beta ^O - \hat{\beta})
   \leq  -\frac{8}{10}c  \lambda\|(\nu)_{S^c}\|_1 + 2\lambda\|(\nu)_{S}\|_1,
\end{equation}
which implies $\lambda\|(\nu)_{S^c}\|_1 \leq \frac{5}{2c}\lambda\|(\nu)_{S}\|_1$ and $\nu \in \mathcal{C}$.
\end{proof}

\subsection{Proof of Theorem \ref{thm:support}}\label{proof:thmsupport}
\begin{proof}
According to Lemma \ref{lem:cvxity}, we will have at $\hat{\beta}$ and $\beta^O$, respectively:
\[f_{\lambda}(\beta^O) \geq f_{\lambda}(\hat{\beta}) + \nabla f_{\lambda}(\hat{\beta})^T(\beta^O - \hat{\beta}) + \frac{\gamma - \eta^-}{2}\|\beta^O - \hat{\beta}\|^2_2,\]
\[f_{\lambda}(\hat{\beta}) \geq f_{\lambda}(\beta^O) + \nabla f_{\lambda}(\beta^O)^T(\hat{\beta} - \beta^O) + \frac{\gamma - \eta^-}{2}\|\hat{\beta} - \beta^O\|^2_2.\]
Since $\ell_1$ norm penalty is convex, we will have
\[\lambda \|\hat{\beta}\|_1 \geq  \lambda\|\beta^O\|_1 + \lambda(\hat{\beta} - \beta^O)^T\xi^O,\]
\[\lambda \|\beta^O\|_1 \geq \lambda \|\hat{\beta}\|_1 + \lambda(\beta^O - \hat{\beta})^T\hat{\xi},\]
where $\hat{\xi}$ and $\xi^O$ are the same as in Lemma \ref{subset}. Combine the above together, we will have:
\[0 \geq  (\nabla f_{\lambda}(\hat{\beta}) + \lambda \hat{\xi})^T(\beta^O - \hat{\beta}) + (\nabla f_{\lambda}(\beta^O) + \lambda \xi^O)^T(\hat{\beta} - \beta^O) + (\gamma - \eta^-)\|\beta^O - \hat{\beta}\|^2_2.\]
Since both $\hat{\beta}$ and $\beta^O$ are d-stationary, we will have
\[(\hat{\beta} - \beta^O)^T \left(\nabla f_{\lambda}(\beta^O) + \lambda \xi^O\right) \geq 0,\]
and $\hat{\beta}$ is d-stationary, there exists a $\hat{\xi} \in \nabla\{\|\hat{\beta}\|_1\}$ such that
\[(\beta ^O- \hat{\beta})^T \left(\nabla f_{\lambda}(\hat{\beta}) + \lambda \hat{\xi}\right) \geq 0.\]
We thus will have $0 \geq (\gamma - \eta^-)\|\beta^O - \hat{\beta}\|^2_2$, which implies that $\beta^O = \hat{\beta}$.
\end{proof}


\section{Proofs in Section \ref{GLOSS}}\label{appB}
\subsection{Proof of Lemma \ref{thm:gloss}}\label{proof:thmgloss}
\begin{proof}
Since $\hat{\beta}$ is a d-stationary solution to Problem (\ref{eq:objGLM}), we have
\begin{equation}
\label{eq:FOC1}
 \nabla L(\hat{\beta}) +  \lambda z -  \nabla h_\lambda(\hat{\beta}) = 0,
\end{equation}
where $z \in \partial g(\hat{\beta})$, where  $\partial g(\hat{\beta})$ is the set of subgradient of $g(\beta) = \|\beta\|_1$ at $\hat{\beta}$.
We can get the gradient for the loss function
$$\nabla L(\hat{\beta}) = \frac{1}{n} \sum_{i = 1}^n (\psi'(X_i^T\hat{\beta})X_i - Y_iX_i).$$
We can further write the above expression as
$$\nabla L(\hat{\beta}) = \frac{1}{n} \sum_{i = 1}^n \big((\psi'(X_i^T\hat{\beta})X_i - \psi'(X_i^T\beta^\ast)X_i) + (\psi'(X_i^T\beta^\ast)X_i - Y_iX_i)\big),$$
where the term $(\psi'(X_i^T\beta^\ast) - Y_iX_i)$ does not depend on the d-stationary solution.
Multiply both side of (\ref{eq:FOC1}) by $(\beta^\ast - \hat{\beta})^T$, we have
\begin{equation}
(\beta^\ast - \hat{\beta})^T\big\{ \frac{1}{n} \sum_{i = 1}^n \big((\psi'(X_i^T\hat{\beta})X_i - \psi'(X_i^T\beta^\ast)X_i) + (\psi'(X_i^T\beta^\ast)X_i - Y_iX_i)\big) +  \lambda z -  \nabla h_\lambda(\hat{\beta}) \big\} = 0.
\end{equation}
By rearranging the terms, we have
\begin{equation*}
\begin{split}
0 &\leq (\beta^\ast - \hat{\beta})^T\big\{ \frac{1}{n} \sum_{i = 1}^n (\psi'(X_i^T\beta^\ast)X_i - \psi'(X_i^T\hat{\beta})X_i)\big\} \\
&= (\beta^\ast - \hat{\beta})^T\big\{\frac{1}{n} \sum_{i = 1}^n(\psi'(X_i^T\beta^\ast)X_i - Y_iX_i)+  \lambda z -  \nabla h_\lambda(\hat{\beta}) \big\}\\
&\leq (\beta^\ast - \hat{\beta})^T\big\{\frac{1}{n} \sum_{i = 1}^n(\psi'(X_i^T\beta^\ast)X_i - Y_iX_i)\big\} + 2\lambda \|(\beta^\ast - \hat{\beta})_S\|_1 \\
&- \lambda \|\hat{\beta}_{S^c}\|_1 + \hat{\beta}_{S^c}^T\nabla h_\lambda(\hat{\beta}_{S^c})\\
&\leq (2+\frac{c}{2})\lambda \|(\beta^\ast - \hat{\beta})_S\|_1 - \frac{c}{2}\lambda \|\hat{\beta}_{S^c}\|_1,
\end{split}
\end{equation*}
where the first ``$\leq$" is due to the convexity of the cumulant function, and the last one is due to the assumptions. We thus conclude that $\hat{\beta} \in \mathcal{C} $, which is defined in the Assumption \ref{ass:RSCGLM}.
Given the restricted strong convexity, according to Lemma \ref{lem:cvxity}, let $f(\beta) = L(\beta) - h_\lambda(\beta)$, we will have
$$f_{\lambda}(\beta^\ast) \geq f_{\lambda}(\hat{\beta}) + \nabla f_{\lambda}(\hat{\beta})^T(\beta^\ast - \hat{\beta}) + \frac{\gamma - \eta^-}{2}\|\beta^\ast - \hat{\beta}\|^2_2$$
and
$$f_{\lambda}(\hat{\beta}) \geq f_{\lambda}(\beta^\ast) + \nabla f_{\lambda}(\beta^\ast)^T(\hat{\beta} - \beta^\ast) + \frac{\gamma - \eta^-}{2}\|\beta^\ast - \hat{\beta}\|^2_2.$$
Adding the above up, we have
\begin{equation}
\label{eq:p2}
\nabla f_{\lambda}(\hat{\beta})^T(\hat{\beta} - \beta^\ast)  \geq  \nabla f_{\lambda}(\beta^\ast)^T(\hat{\beta} - \beta^\ast) + (\gamma - \eta^-)\|\beta^\ast - \hat{\beta}\|^2_2.
\end{equation}
Adding $\lambda z^T(\hat{\beta} - \beta^\ast)$ to both side, we will have
\begin{equation}
\label{eq:p1}
0 = (\nabla f_{\lambda} (\hat{\beta})^T+ \lambda z^T)(\hat{\beta} - \beta^\ast)  \geq  \nabla f_{\lambda}(\beta^\ast)^T(\hat{\beta} - \beta^\ast) + \lambda z^T(\hat{\beta} - \beta^\ast) + (\gamma - \eta^-)\|\beta^\ast - \hat{\beta}\|^2_2,
\end{equation}
From inequalities (\ref{eq:p1}) and (\ref{eq:p2}), we have
\begin{equation}
\begin{split}
(\gamma - \eta^-)\|\beta^\ast - \hat{\beta}\|^2_2
&\leq -\nabla f_{\lambda}(\beta^\ast)^T(\hat{\beta} - \beta^\ast) - \lambda z^T(\hat{\beta} - \beta^\ast)\\
&\leq (2+\frac{c}{2})\lambda \|(\beta^\ast - \hat{\beta})_S\|_1 - \frac{c}{2}\lambda \|\hat{\beta}_{S^c}\|_1\\
&\leq (2+\frac{c}{2})\lambda\sqrt{|S|}\|\hat{\beta} - \beta^\ast\|_2
\end{split}
\end{equation}
where the last ``$\leq$" is due to the fact that $\|(\hat{\beta} - \beta^\ast)_S\|_1 \leq \sqrt{|S|}\|\hat{\beta} - \beta^\ast\|_2$.
We thus derive the bound that
\begin{equation}
\begin{split}
\|\beta^\ast - \hat{\beta}\|_2 \leq \frac{(4+c)\lambda}{2(\gamma - \eta^-)}\sqrt{|S|}
\end{split}
\end{equation}
\end{proof}

\section{Proofs in Section \ref{SEC:CCCP}}\label{appC}
\subsection{Proof of Lemma \ref{lemma:Decrease}}\label{proof:lemmaDecrease}
\begin{proof}
Given $\beta_{k}$ as the update in the kth iteration, we adopt the following procedure to update the estimation:
\begin{equation}
\label{eq:updateDecrease}
\beta_{k+1} = \arg \min L_n(\beta) + \sum_{i = 1}^p (\lambda - h'(|\beta_{ki}|)) |\beta_i|.
\end{equation}
Let $Q(\beta|\beta_{k}) = L_n(\beta) + \lambda\|\beta_{k}\|_1 - h_{\lambda}(\beta_{k}) + \sum_{i = 1}^p (\lambda - h'(|\beta_{ki}|)) (|\beta_i| - |\beta_{ki}|)$. It can be easily checked that  $Q(\beta_k|\beta_{k}) = F(\beta_k)$ and the following is equivalent to \ref{eq:updateDecrease}:
\begin{equation}
\beta_{k+1}  = \min Q(\beta|\beta_{k}).
\end{equation}
Since $h_{\lambda}(\cdot)$ is convex by Assumption \ref{ass:hcvxity}, we have
$$F(\beta) \leq Q(\beta|\beta_{k}).$$
According to the definition of $\beta_{k+1}$, we have
$$F(\beta_{k+1}) \leq Q(\beta_{k+1}|\beta_{k}) \leq Q(\beta_{k}|\beta_{k}) = F(\beta_k).$$
\end{proof}

\end{document}